\documentclass{article}

\usepackage{arxiv}

\usepackage[utf8]{inputenc} 
\usepackage[T1]{fontenc}    
\usepackage[hidelinks]{hyperref}       
\usepackage{url}            
\usepackage{booktabs}       
\usepackage{amsfonts}       
\usepackage{nicefrac}       
\usepackage{microtype}      
\usepackage{xcolor, colortbl}   
\usepackage{amssymb, amsmath}
\usepackage{float}
\usepackage{enumitem}
\usepackage{mathtools, cuted}
\usepackage{verbatim}
\usepackage{tikz,xcolor}
\usepackage{stfloats}
\usepackage{sidecap}
\usepackage{algorithm}
\usepackage[italicComments=false, commentColor=darkgray]{algpseudocodex}   
\usepackage{array,ragged2e,multirow}
\usepackage{scrtime}
\usepackage[usestackEOL]{stackengine}
\usepackage[bb=dsserif]{mathalpha}
\usepackage[title]{appendix}

\tikzset{algpxIndentLine/.style={draw=lightgray,very thin}}  
\usetikzlibrary{shapes.geometric, arrows, positioning}

\definecolor{lightblue}{HTML}{DAE3F3}  
\definecolor{darkblue}{HTML}{41719C}

\newcommand{\pinP}{p\in \mathcal{P}}
\newcommand{\dinD}{d \in \mathcal{D}}
\newcommand{\finF}{f \in \mathcal{F}}
\newcommand{\winW}{w \in \mathcal{W}}
\newcommand{\minM}{m \in \mathcal{M}}

\newcommand{\iinK}{i \in \mathcal{K}_p}
\newcommand{\R}{\mathbb{R}}

\newcommand{\D}{\mathcal{D}}
\newcommand{\Bm}{\mathcal{B_M}}
\newcommand{\bm}{b_\mathcal{M}}
\newcommand{\binBm}{b_\mathcal{M} \in \mathcal{B_M}}

\newcommand{\Mh}{\mathcal{M}_h}
\newcommand{\Mhpref}{\mathcal{M}_h^{\text{pref}}}
\newcommand{\Mhallow}{\mathcal{M}_h^{\text{allowed}}}
\newcommand{\Ppref}{\mathcal{P}^{\mathit{w{\text{ - pref}}}}}
\newcommand{\pinPpref}{p\in \Ppref}
\newcommand{\wPref}{w^{\text{pref},p}}
\newcommand{\minbM}{m \in b_\mathcal{M}}
\newcommand{\fminH}{\mathit{f}_{h}^{\min}}
\newcommand{\fminHp}{\mathit{f}_{h_p}^{\min}}
\newcommand{\ginG}{g \in \mathcal{G}}
\newcommand{\Gubm}{\mathcal{G}^\mathcal{U}_{b_\mathcal{M}}}

\newcommand{\dMin}{d_{p}^{\min}}
\newcommand{\SMpref}{\mathcal{S_M}^{\text{pref,}p} }
\newcommand{\Hlist}{\mathcal{H}_{\mathit{list}}}
\newcommand{\Plist}{\mathcal{P}^h_{\textit{list}}}
\newcommand{\Mlist}{\mathcal{M}^p_{\textit{list}}}
\newcommand{\Wlist}{\mathcal{W}^p_{\textit{list}}}
\newcommand{\StempMDW}{S^{\textit{temp}}_{m,d,w}}
\newcommand{\Stemp}{S^{\textit{temp}}}
\newcommand{\dstart}[1]{d^i_{\text{start}, #1}}
\newcommand{\dend}[1]{d^i_{\text{end}, #1}}
\newcommand{\Ahp}{\mathcal{A}_{h_p}}
\newcommand{\DmU}{\mathcal{D}_m^{\mathcal{U}}}
\newcommand{\DmUhat}{\mathcal{D}_{\hat{m}}^{\mathcal{U}}}
\newcommand{\durpfirst}{\mathit{dur}_{p}^{1}}
\newcommand{\durprest}{\mathit{dur}_{p}^{f}}
\newcommand{\q}[2]{q^i_{p,m,#1,#2}}

\newcolumntype{P}[1]{>{\RaggedRight\arraybackslash}p{#1}}
\newcommand{\mypar}[1]{\paragraph{#1.}}

\renewcommand{\headeright}{}
\renewcommand{\undertitle}{}
\renewcommand{\shorttitle}{}

\tikzstyle{node} = 
[rectangle, rounded corners, 
minimum width=2.8cm, 
minimum height=3.8cm,
text centered, 
text width=2.8cm, 
draw=darkblue,line width=0.3mm, 
fill=lightblue,
font={\fontsize{9pt}{10}\selectfont}]  

\tikzstyle{widenode} = 
[rectangle, rounded corners, 
minimum width=3.1cm, 
minimum height=3.8cm,
text centered, 
text width=3.1cm, 
draw=darkblue, line width=0.3mm, 
fill=lightblue,
font={\fontsize{9pt}{10}\selectfont}]

\tikzstyle{thinnode} = 
[rectangle, rounded corners, 
minimum width=2.0cm, 
minimum height=3.8cm,
text centered, 
text width=2.0cm, 
draw=darkblue, line width=0.3mm, 
fill=lightblue,
font={\fontsize{9pt}{10}\selectfont}]

\tikzstyle{smallnode} = 
[rectangle, rounded corners, 
minimum width=5cm, 
minimum height=0.8cm,
text centered, 
text width=5cm, 
draw=darkblue, line width=0.3mm, 
fill=lightblue,
font={\fontsize{9pt}{10}\selectfont}]

\tikzstyle{arrow} = [line width=0.3mm,->,>=stealth, darkblue, text=black]

\title{A Column Generation Approach for Radiation Therapy Patient Scheduling with Planned Machine Unavailability and Uncertain Future Arrivals}

\author{
  Sara Frimodig\\
  RaySearch Laboratories\\
  Department of Mathematics \\
  KTH Royal Institute of Technology\\
  Stockholm, Sweden \\
  \texttt{sarhal@kth.se} \\
  \And
  Per Enqvist \\
  Department of Mathematics \\
  KTH Royal Institute of Technology\\
  Stockholm, Sweden \\
  \And
  Jan Kronqvist \\
  Department of Mathematics \\
  KTH Royal Institute of Technology\\
  Stockholm, Sweden \\
}

\begin{document}
\maketitle

\begin{abstract}
The number of cancer cases per year is rapidly increasing worldwide. In radiation therapy (RT), radiation from linear accelerators is used to kill malignant tumor cells.  Scheduling patients for RT is difficult both due to the numerous medical and technical constraints, 
and because of the stochastic inflow of patients with different urgency levels. 
In this paper, a Column Generation (CG) approach is proposed for the RT patient scheduling problem. The model includes all the constraints necessary for the generated schedules to work in practice, including for example different machine compatibilities, individualized patient protocols, and multiple hospital sites. The model is the first to include planned interruptions in treatments due to maintenance on machines, which is an important aspect when scheduling patients in practice, as it can create bottlenecks in the patient flow. Different methods to ensure that there are available resources for high priority patients at arrival are compared, including static and dynamic time reservation. Data from Iridium Netwerk, the largest cancer center in Belgium, is used to evaluate the CG approach. The results show that the dynamic time reservation method outperforms the other methods used to handle uncertainty in future urgent patients. A sensitivity analysis also shows that the dynamic time reservation method is robust to fluctuations in arrival rates. The CG approach produces schedules that fulfill all the medical and technical constraints posed at Iridium Netwerk with acceptable computation times.

\end{abstract}


\section{Introduction}
Cancer is one of the leading causes of premature mortality worldwide. By 2040, the predicted number of new cancer cases per year is expected to exceed 27 million \cite{Wild2020}, a 40\% increase compared to the estimated 19.3 million cancer cases in 2020 \cite{Ferlay2021}. Radiation therapy (RT) is a cancer treatment that uses high doses of radiation to destroy or damage cancer cells.
As a consequence of the rising cancer incidents, the need for RT will increase \cite{Borras2016}. 
In external beam radiation, a machine called linear accelerator (\emph{linac}) is used to direct radiation from outside the body into the tumor. 
Each treatment is either curative, with intent to cure the patient, or palliative, with intent to improve quality of life by providing symptom control. Cancer patients are often divided into different urgency levels depending mainly on the site of the cancer and treatment intent. For high priority patients, the treatment should start as soon as possible after admission, while lower prioritized patients can wait for two to four weeks. 

The duration of each treatment session varies between patients due to for example treatment technique and tumor location. RT treatments are normally divided into a number of sessions called \emph{fractions} that are delivered daily over several weeks. Delivering a small fraction of the total radiation dose allows time for normal cells to repair between the treatments, which reduces the side effects. However, the overall treatment time should be kept as short as possible, as longer gaps between fractions can enable the repopulation of cancer cells to accelerate, leading to potentially lower cure rates \cite{RCR2019}. Patients usually receive treatment five days a week, but if there is some machine unavailability it can be necessary to postpone fractions, leading to gaps in the schedule.

In the RT scheduling problem (RTSP), the aim is to schedule patients for RT, given a set of linacs, for a certain planning horizon. 
Patients can have preferences on what time during the day they want to be treated, and on what hospital for cancer centers with multiple sites. Moreover, the treatments have different machine requirements. Other difficulties include patients that are treated with multiple consecutive treatments, or with non-conventional treatments. Furthermore, planned and unplanned unavailability of the machines can create bottlenecks in the patient flow.
One of the main challenges for the RTSP is that there is uncertainty in demand. Since the patients are of different priority, it is important that there are resources available for urgent patients at arrival. In practice, this is often handled by reserving a percentage of the machine capacity for high priority patients. This method can cause delays in treatments, as well as unnecessary idle time on the machines, especially at large clinics with high patient flow. 

Long waiting times for RT negatively impacts clinical outcomes. For example, it can cause patients to experience prolonged symptoms, tumor growth, and psychological distress  \cite{Zumer2020,CHEN2008,Fortin2002,Gomez2015,VanHarten2015}.  Long waiting times has also been identified to cause stress for RT staff, which can compromise quality and safety of the treatments \cite{French2004}. The waiting time for RT is often directly linked to the RTSP, since the number of linacs in a clinic is usually limited.  At cancer clinics today, almost all patient scheduling is done manually by the staff. 
As the demand for RT grows, efficient resource planning is an important tool to achieve short waiting times.  Therefore, this paper makes the following contributions:
\begin{itemize}
\item The \emph{main contribution} is that we present an automatic scheduling algorithm for the RTSP, with all constraints and objectives necessary for the model to work in practice. For the first time, planned machine unavailability is included in an RTSP model, which is an important step towards a full clinical implementation. We also present a method for dynamic time reservation to handle uncertainty in future patient arrivals. 
\item The method is evaluated using data from Iridium Netwerk, a large RT center located in Antwerp, Belgium. In 2020, they operated 10 linacs, delivering 5500 RT treatments to approximately 4000 patients. 
\item The main \emph{technical novelty} lies in the column generation (CG) approach. To the best of our knowledge, this is the first model to simultaneously assign all fractions of the patients to both linacs and specific time windows, while also considering \textit{planned machine unavailability} and the constraints and objectives related to the resulting gaps in the schedules. It is also the first model to include \textit{consecutive treatments}, where a primary treatment is followed by a secondary treatment with some additional constraints, and \textit{non-conventional treatments}, such as treatments that should be delivered every-other-day. The model also supports \textit{multiple hospital locations} and allows the patients to have hospital site preferences. Furthermore, the model includes all the medical and technical constraints necessary for the scheduling to work in practice at Iridium Netwerk. 
\end{itemize}
The paper is organized as follows. Section~\ref{Sec:rel_work} presents related work, followed by a problem description in Section~\ref{Sec: Problem Formulation}. Section~\ref{Sec:ColGen} presents the column generation model. Methods to manage uncertainty in patient arrivals are discussed in Section~\ref{Sec:FutureArrivals}.  Section~\ref{Sec:results} presents  the experiments and the numerical results. The conclusions are presented in Section~\ref{Sec:conclusions}. 

\section{Related Work}
\label{Sec:rel_work}
Scheduling in healthcare has been widely studied, and summarized in several extensive literature reviews. Cayirli et al. \cite{Cayirli2003}, Gupta et al. \cite{Gupta2008}, and Ahmadi-Javid et al. \cite{Ahmadi-Javid2017} present reviews that focus on scheduling on a single resource, and Marynissen et al.  \cite{Marynissen2019} review the special case of multi-appointment scheduling. 
Patient scheduling consists of \emph{allocation scheduling} and \emph{appointment scheduling}. Allocation scheduling refers to methodologies for allocating patients to resources in advance of the service date, when future demand is still unknown, without assigning specific appointment times.  In contrast, in appointment scheduling all patients for a given service day are assumed to be known, and specific resources and starting times are assigned to the patients. Appointment scheduling problems 
for example aim to minimize machine idle time, maximize preference satisfaction regarding treatment time, or deal with uncertain treatment durations or delays. On the other hand, allocation scheduling must fulfill capacity constraints, and usually deals with patients of different types and priorities. Most allocation scheduling algorithms are intended to dynamically allocate resources to patients using a rolling time horizon. The problem studied in this paper aims to bridge the gap between appointment and allocation scheduling by performing these two scheduling tasks simultaneously.





The RTSP has been studied in various alternations in the past 15 years. For a review of the literature, see 
Vieira~et~al.~\cite{Vieira2016}, in which the authors found $12$~papers addressing the problem of scheduling RT patients on linacs. The solution methods have varied; exact methods such as integer programming (IP) have been used to solve small instances and metaheuristics have been used for larger ones. 

Patient scheduling is done either in an \emph{online} fashion, where each request is handled immediately, or in \emph{batches}, where the scheduling is done at certain time intervals (such as daily or weekly). Using batch scheduling, Conforti et al. \cite{Conforti2008} present the first IP model for optimization of RT appointments using treatment slots of equal length (\emph{block} scheduling), a model that was later developed by the same authors to allow for different fraction durations (\emph{non-block} scheduling) \cite{Conforti2010}.  Jacquemin et al. \cite{Jacquemin2011} introduce the notion of treatment patterns in an IP model to allow non-consecutive treatment days. These papers all present models for a short planning horizon that do not consider all the constraints present in real-world RT scheduling, such as multiple machine types and partial availability in the schedule. For a more realistic setup, Petrovic et al. \cite{Petrovic2006,Petrovic2013} propose heuristic and metaheuristic approaches for block scheduling based on prioritized rules, where the latter considers both scheduling of the pretreatment phase and the treatment phase. 

The first method for dynamic RT scheduling taking future events into account was presented by Sauré et al. \cite{Saure2012}, where the problem is modeled as a discounted infinite-horizon Markov decision process that finds an approximately optimal scheduling policy. Their proposed policy can increase the  treatments initiated within 10 days from 73\% to 96\% compared to a myopic policy (i.e., not taking future patients into account). Gocgun \cite{Gocgun2018} later extended the same problem setup to also include patient cancellations. These papers assume a simplified model of a cancer center, equipped with three identical machines, 8.25 requests per day and scheduling done in batches. In contrast, Legrain et al. \cite{Legrain2015} propose a hybrid method combining stochastic and online optimization to dynamically schedule patients as they arrive. Also their setup is small; they consider block scheduling with two linacs and less than 3.5 requests per day. Aringhieri et al. \cite{Aringhieri2020} also present methods for online RT scheduling, and develop three online optimization algorithms for a block-scheduling formulation and one machine. 
These methods all allocate a start day to each patient, with no sequencing of patients throughout the day, and are difficult to scale to realistic instances.


Particle therapy (PT) is a form of RT where a single particle beam is shared between multiple treatment rooms. The medical and technical constraints differ from conventional photon beam RT, and most studies focus on maximizing the beam utilization by optimal appointment scheduling. Maschler et al. \cite{Maschler2016} showed that the exact formulation of the problem is highly intractable. Using different heuristic methods, Vogl et al. \cite{Vogl2019} and Maschler et al. \cite{Maschler2020} both create a schedule with treatments close to a pre-defined target time. Accounting for uncertain activity durations, Braune et al. \cite{Braune2021} present a stochastic optimization model for appointment scheduling in PT and solve it using heuristics.

Focusing on appointment scheduling, where the patient list is assumed to be known and the main task is the sequencing of patients throughout the day, Vieira et al. \cite{Vieira2020}  create weekly schedules using a mixed-integer programming (MIP) model combined with a pre-processing heuristic that divides the problem into subproblems for clusters of machines. Their objective is to maximize the fulfillment of the patients' time window preferences. In \cite{Vieira2021}, the same authors evaluate their method in two Dutch clinics, showing that the weekly schedule was improved in both centers. 
Another paper that focuses on patient sequencing is Moradi et al. \cite{Moradi2022}, where the authors present a data-driven approach that uses the patient information to improve the weekly schedules. The predictions are utilized in a MIP model to determine the optimal sequence of patients, for a list of patients that have previously been assigned a day. The model is simplified; all patient durations are equal and all machines are identical and independent. However, the results seem promising; it is favorable to schedule reliable patients early on to reduce idle time on machines caused by delayed patients or no-shows.

A two-stage approach for the RTSP is proposed by Pham et al. \cite{Pham2021}. In the first phase, an IP model assigns patients to linacs and days, and the second phase assigns specific appointment times using either a MIP or a constraint programming (CP) model. 
The authors evaluate the algorithm dynamically on a rolling time horizon using different scheduling strategies. The approach does not take future arrivals into account when making scheduling decisions; a certain percent of linac capacity is saved for urgent patients in the first phase. 
Their approach is tested using generated data with seven linacs and a time horizon of 60 days based on data from CHUM, a cancer center in Canada. They show that the CP model finds good solutions sooner, while the MIP model proves optimality faster. 

In Frimodig et al. \cite{Frimodig2022}, the sequencing of patients is decided at the same time as the assignment to linacs and days. The authors compare an IP model, a CP model and a column generation (CG) model to solve the problem, and include the expected future patients to dynamically reserve linac capacity for future urgent patients. The models are tested on generated data based on data from Iridium Netwerk in Belgium, that has ten linacs. The results show that the CG model outperforms the other in all problem instances. The setup is similar to the setup in this paper, with the difference that the models do not consider unavailability of machines, consecutive treatments, or non-conventional treatments. Furthermore, they do not evaluate the CG approach dynamically on a rolling time horizon; that is done in this paper.



CG is a decomposition technique that if often successfully used for solving huge integer programs. The method was first presented by Ford and Fulkerson \cite{FordFulkerson1958}, and has the advantage that it does not consider all variables explicitly, but instead only generates the variables that have the potential to improve the objective function. The method alternates between a restricted master problem (RMP) and one or more subproblems used to generate new variables to the RMP. 
When applying CG in scheduling problems, the decision variable in the master problem is most often a schedule for one day, or one shift, or one person, and the master problem is a set partitioning problem used to find the optimal overall schedule by minimizing the sum of the costs of the associated variables. The subproblems are used to find the variables. Other medical fields where CG has been applied is for surgeon and surgery scheduling \cite{Wang2018}, for patient admission \cite{Range2014}, and for nurse scheduling \cite{Bard2005}. In cancer treatments, CG has also been used for brachytherapy scheduling using deteriorating treatment times \cite{Shao2021}, and for intraday scheduling in chemotherapy \cite{Lyon2022}.

\section{Problem Description}
\label{Sec: Problem Formulation}
\begin{figure*}[ht]
    \center
    \includegraphics[width=16.5cm]{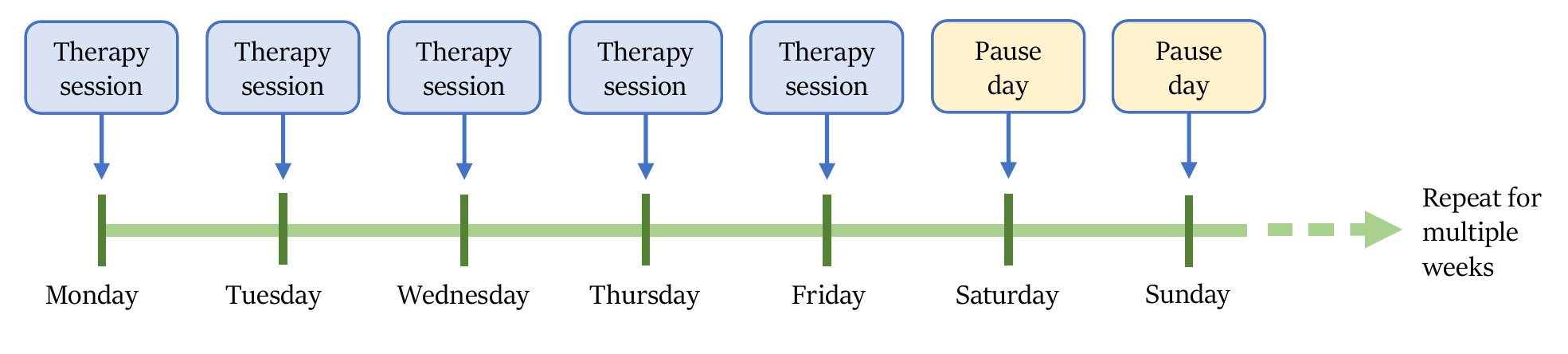}
    \vspace{-0.8cm}
    \caption{A typical fractionation scheme of an RT patient}
    \label{fig:treatment_course}
\end{figure*}

The task in the RTSP is to assign patients to treatment machines (linacs) for each fraction according to a specific set of rules and objectives. This section presents the real-world constraints and objectives present at Iridium Netwerk. Figure~\ref{fig:treatment_course} presents how a patient usually is treated with RT.

The instructions for how a patient should be scheduled are given primarily by the \emph{treatment protocol}, which is assigned to the patient by the treating physician. The most common treatment protocols at Iridium Netwerk can be seen in Table \ref{tab:protocols}. Each treatment protocol is associated with a \emph{priority} based on the urgency for treatment and the treatment intent. There are three priority groups: A, B and C. In 2020 at Iridium, approximately $37\%$ of the patients were priority~A, $16\%$ were priority~B, and $46\%$ were priority~C. The treatment protocol states the \emph{minimum number of fractions per week}, that is, how much it is allowed to deviate from the ideal five fractions per week.  The protocol also states the \emph{preferred machines} for the treatment, and the machines that are \emph{allowed} (but not preferred). Furthermore, the protocol states the minimum number of \emph{days from CT}, which is the time period from the mandatory CT simulation to start of RT treatment used to create the treatment plan. The protocols have \emph{allowed start days}: palliative patients can start any weekday, whereas curative patients cannot start on Fridays. In addition to the treatment protocol, there is also patient specific information: the \emph{number of fractions} during which the patient will be treated, the \emph{dose prescription} for each fraction, and the \emph{duration} of the first and subsequent fractions, where the first is longer than the rest for each patient because of initial setup and quality checks, as well as extra time for patient education and reassurance. 

Iridium Netwerk is a network with four different hospitals that each have between two and four linacs. Two linacs are \emph{completely beam-matched} if they are the same machine type at the same hospital, and \emph{partially beam-matched} if they are the same machine type, but at different hospitals. Switching between completely beam-matched machines between fractions can be done at no cost, whereas there is a cost for switching to a machine that is only a partially matched. 
\begin{table}[ht]
\caption{The most common treatment protocols Iridium Netwerk}
\label{tab:protocols}
\begin{tabular}{P{1.7cm} P{1.05cm} P{2.1cm}  P{1.97cm} P{4.5cm} P{2.66cm}}
\hline\noalign{\smallskip}
Protocol & Priority & Minimum fractions/week & Minimum days from CT & Preferred machines &Allowed machines\\
\hline\noalign{\smallskip}
Palliative & A & 1 & 0 & M1, M2, M3, M4, M5, M6, M7, M8 & M10\\
Breast & C & 3 & 7 & M1, M4, M5, M6, M8 & M2, M3, M7, M10 \\
Prostate & C & 3  & 9 & M1, M3, M4, M5, M6, M7, M8 & M2, M10 \\
Head-Neck  & A &  5 & 11 & M2, M3, M5, M6, M10 & M1, M4 \\
Lung & B & 4  & 9  & M2, M3, M5, M6, M7, M10 & M1, M4, M8 \\
\hline\noalign{\smallskip}
\end{tabular}
\vspace{-8pt}
\end{table}


The day is divided into four \emph{time windows},
and the model assigns each fraction to a machine, time window and day. 
Assigning patients to time windows instead of exact starting times leads to a more efficient model, while it enables all the objectives and constraints needed to follow the instructions from the scheduling staff at Iridium Netwerk.

Priority B and C patients are assumed to be notified of their starting day one week (five working days) in advance, which the majority of patients find reasonable according to the literature \cite{Olivotto2015}. Priority~A patients are notified immediately. All fractions are communicated at once, as this is the current practice at Iridium. Each patient's schedule can change until is has been communicated, i.e., booking decisions are postponed to the next day if patients are scheduled after the notification period. A \emph{daily batch scheduling} strategy is applied, using information accumulated throughout the day.

\begin{figure}
\centering
\includegraphics[width=0.99\linewidth]{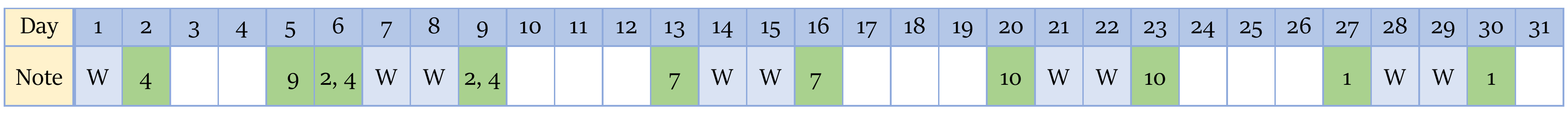}
\vspace{-3pt}
\caption{Schedule for March. "W" indicates weekend and a number indicates planned unavailability on that machine}
\label{Fig:MarchSchedule}
\vspace{-0.4cm}
\end{figure}

For an automatic scheduling algorithm to work in practice, machine unavailability is something that needs to be considered in some way. In \cite{RCR2019}, the UK guidelines for the maximum overall treatment times are presented. The gaps between fractions are often caused by machine failure or some other unexpected event, and the gaps induced by the scheduling should be kept to a minimum. Figure~\ref{Fig:MarchSchedule} shows when machines have \emph{planned unavailability} for an example month in 2020 at Iridium Netwerk. Depending on each patient's treatment protocol, 
the affected fractions can either be scheduled on a beam-matched machine the same day, or postponed by adding them after the last one, or re-scheduled the same week by scheduling two fractions on some day the same week as the unavailability. 

Some patients are treated with two or three \emph{consecutive treatments}, where the primary treatment is directly followed by a secondary treatment. For example, breast cancer patients are often treated with a boost plan that follows after the primary plan. The secondary treatment must be handled separately in the scheduling since both machine requirements and durations of the fractions can differ from the primary plan. The secondary plan also has the extra constraints that the first day should follow after the last of the primary plan (ideally start the day after the primary plan has ended). At Iridium Netwerk in 2020, around 17\% of the patients were treated with two or more consecutive treatments. 
Furthermore, there are so called \emph{non-conventional treatments}, meaning treatments that do not follow the regular five fractions per week schedule. The most common are treatments that should be delivered three times per week, with a pause day between each fraction. At Iridium in 2020, around 9\% of all treatments were non-conventional.

Because of  \emph{uncertainty in the arrival rates}, it is important that there are resources available for urgent patients at arrival. This is difficult since the treatments span multiple weeks (often $5-7$ weeks) and the patients have different priority. At Iridium Netwerk, this is handled by the booking administrator, who reserves empty timeslots for urgent patients on each machine. In this paper, different methods for handling uncertainty in future arrivals are presented in Section~\ref{Sec:FutureArrivals}.

The \emph{scheduling objectives} are formulated in collaboration with Iridium Netwerk. The most important objective is to minimize the waiting times, especially for urgent patients. The patients should ideally be scheduled around the same time every day. For patients that have a preference for treatment time, it should try to be fulfilled. The treatment should if possible be scheduled on a preferred machine, and the number of switches to a machine that is not beam-matched should be kept to a minimum. Since Iridium has multiple hospital sites, the schedule should try to meet the patient preference on treating hospital. Finally, the overall treatment time for each patient should be kept as short as possible to avoid gaps in the schedule.  The objectives are combined into an objective function that is presented in Section~\ref{Sec:subproblems}. 


\section{Column Generation Model}
\label{Sec:ColGen}
\begin{figure}[ht]
\begin{tikzpicture}[scale=1.0,transform shape, align=center]
\node (start) [node] {Initialize restricted master problem (RMP) by generating a reduced set of schedules for each patient~$p$ (Section~\ref{Sec:CG-heuristic})};

\node (master) [node, right=0.6cm of start] {\emph{Master problem: \\Schedule selection} \\
Solve LP relaxation of RMP to select schedules and generate dual variables
(Section~\ref{Sec:CG_master})};

\node (sub) [widenode, right=0.6cm of master] {\emph{Subproblem: \\ Schedule generation} \\
Update subproblems, one for each patient, with dual variables and solve resulting IPs to generate new columns (schedules)\\
(Section~\ref{Sec:subproblems})};  

\node (cost) [thinnode, right=0.6cm of sub] {Are any reduced cost negative?};

\node (addcolumn) [smallnode, below=0.5 cm of sub] {Add new columns to RMP};

\node (stop) [thinnode, right=0.6cm of cost] {Solve most recent RMP using integer values to select final schedule};

\draw [arrow] (start) -- (master);
\draw [arrow] (master) -- (sub);
\draw [arrow] (sub) -- (cost);
\draw [arrow] (cost) -- node[anchor=south] {No} (stop);
\draw [arrow] (cost) |- node[anchor=west] {Yes} (addcolumn);
\draw [arrow] (addcolumn) -| (master);
\end{tikzpicture}
\caption{The column generation algorithm}
\label{Fig:ColGenFlow}
\end{figure}
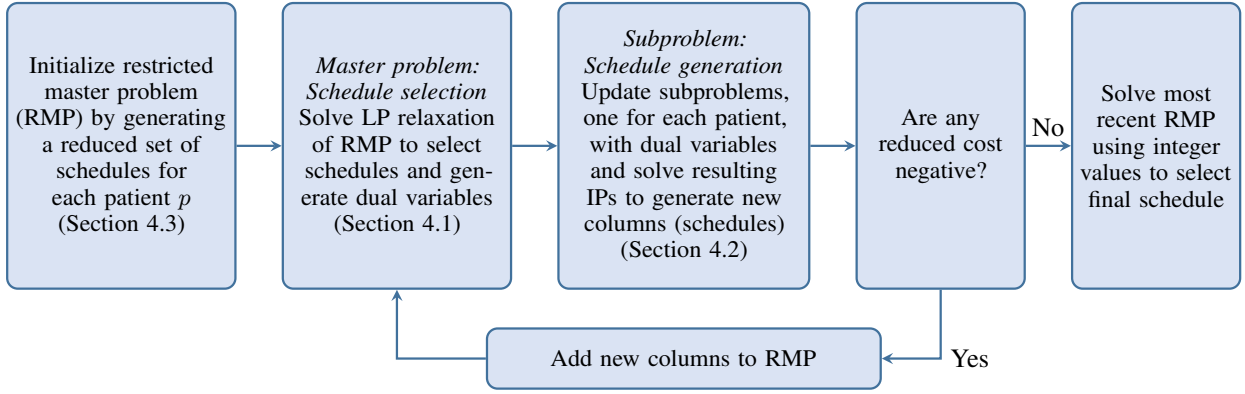

The RTSP is formulated as a binary set partitioning model, where each decision variable is a schedule for a patient, and the aim is to choose the optimal set of patient schedules. Since it would be too expensive to generate all feasible patient schedules, a CG approach is used, which consists of a (restricted) master problem and one subproblem for each patient. The master problem is used for \emph{schedule selection}: it is solved to choose a schedule for each patient to make the overall schedule feasible and optimal. The subproblems are used for \emph{schedule generation}: for each patient, a new schedule is generated that fulfills all medical and technical constraints, and if the schedule has a negative reduced cost, the variable is added to the master problem.  The CG algorithm is presented in  Figure~\ref{Fig:ColGenFlow}. The notations are introduced in Table~\ref{tab:inputs}.  

The subproblems are isolated from each other and are often very easy to solve. However, the optimization is exclusively guided by the values of the dual variables, which might lead to a large number of iterations in the CG algorithm. To improve the speed of the CG algorithm, all subproblems are run in the first iteration, and henceforth in every third iteration. In the two intermediate iterations, only the subproblems that previously have had a negative reduced cost are run, since these are the most likely to have a negative reduced cost again. 
The CG algorithm terminates when no negative reduced cost variables are generated by the subproblems, which means that the linear relaxation to the original problem has obtained an optimal solution. To get an integer solution the CG procedure is normally embedded in a branch-and-price algorithm (see e.g. \cite{Barnhart1998}), but  this is not done in this paper due to limitations in solution time. Instead, the linear program is converted to an integer program in the last step. This does not ensure an optimal solution, since some schedules not generated by the procedure could potentially improve the integer solutions. 


\begin{table}
\caption{Notations for the models}
\label{tab:inputs}
\begin{small}
\begin{tabular}{P{0.36\linewidth} P{0.58\linewidth}}
\hline\noalign{\smallskip}
    Parameter                                  & Description \\
    \hline\noalign{\smallskip}
	$\mathcal{P} = \{1, \dots, P\}$               & Set of all patients  \\
	$(p_{c,1},p_{c,2}) \in \mathcal{P}^{con}  $                              & List of tuples of connected patients for consecutive treatment $c$\\
	$\D = \{1, \dots, D \}$                          & Set of weekdays in the planning horizon \\
	$\mathcal{W} = \{1,2, 3, 4\}$                 & Set of time windows in a day  \\
	$L_w = \{ 135, 120, 120, 135\}$            & The window length in minutes for window $\winW$\\
		$\mathcal{H}=\{1,\dots,H\} $            & Set of treatment protocols, where $h_p \in \mathcal{H}$ is the protocol of patient $p$ \\
	$\mathcal{M} = \{M1, \dots, M10\}$   & Set of machines \\
	$\Mh = ( \Mhpref \cup  \Mhallow )	\subseteq \mathcal{M}$                                                  & Set of machines for protocol $h$, where  $\Mhpref$ are the set of preferred machines and $\Mhallow$ are the set of allowed (but not preferred) machines \\
	$\mathcal{C_M} = \{\{M3, M10 \}, \{M5, M6\} \}$    & Sets of completely beam-matched machines\\
	$\Bm =  \{ \{M1, M4, M8 \},\{M9\} ,$                     & Sets of all beam-matched machines\\
	$\quad \quad \quad \enskip \{M2, M3, M5, M6, M7, M10\}\}$ & \\
	$\mathcal{S_M} = \{ \{M3, M4, M9, M10\},  \{M1, M2\},$  & Sets of machines at the different hospital sites $S1, \dots, S4$. $\SMpref$  is the  \\
	$\quad \quad \quad \enskip\{M7, M8\}, \{M5, M6\}\} $ & set of machines at site preference of patient $p$\\
    $\durpfirst, \durprest$                     & Duration of first and subsequent fraction for patient $p$ (minutes) \\
	$\mathcal{F}_p \in \{1,\dots, F_{p} \}$     & Set of all fractions for patient $p$	\\
    $S_{m,d,w} \in \{0,\dots,L_w\}$     & Occupied minutes each window $\winW$, machine $\minM$ and day $\dinD$ \\
	$\Ahp \in \D$                                  & Set of allowed start days for the protocol $h$ of patient $p$ \\
	$c_{h_p} \in \{10,3,1\}$                       & Weights for patient $p$'s protocol $h$ based on priority group A, B or C \\
	$\dMin \in \mathcal{D}_w$             & The earliest day for patient $p$ to be scheduled \\
	$\Ppref \subset \mathcal{P}$           & Set of patients that have a time window preference \\
	$\wPref \in \mathcal{W}$   & The window preference of patient $\pinPpref$ \\
	$\fminH \in \{1,\dots,5\}$ 					& Minimum number of fractions per week for protocol $h$ \\
	$\mathcal{G} = \{1,\dots,G\}$				& Set of all weeks in planning horizon. $\mathcal{D}_g$ is the set of days in week $\ginG$.\\
	$\Gubm \subset \mathcal{G}$              & Set of weeks where there is unavailability for machine group $\binBm$\\
	$\DmU \subset \D$                      & Set of unavailable days for machine $m$ \\
\noalign{\smallskip}\hline
\end{tabular}
\end{small}
\end{table}


\subsection{Master Problem: Schedule Selection}
\label{Sec:CG_master}

 In the master problem, each patient has an associated index set $\mathcal{K}_p$ of feasible schedules, and the variable $a^i_{p}=1$ if schedule $\iinK$ is allocated to $\pinP$, and 0 otherwise. Model \eqref{Eq:master obj}-\eqref{Eq:master integer} is the master problem: the restricted master
problem (RMP)  is made of a subset $\mathcal{K}_p' \subset \mathcal{K}_p$ of feasible schedules for each $\pinP$. Each schedule $i \in \mathcal{K}_p $  has a cost $c^i_{p}$, a treatment duration $D^i_{p,m,d,w}$  for each machine, day, and time window, and a start day $\dstart{p}$ and end day $\dend{p}$ representing the first and last days of treatment, computed from the subproblem variables in \eqref{Eq:colgen sub cost}, \eqref{Eq:Dur_master}, \eqref{Eq:dstart} and \eqref{Eq:dend}. 

The objective function \eqref{Eq:master obj} is to minimize the total cost of the chosen schedules. Constraint \eqref{Eq:master a=1} states that exactly one schedule is chosen for each patient. Constraint \eqref{Eq:master sched} ensures that all chosen schedules will fit in the overall schedule. To handle \emph{consecutive treatments}, where a patient has two treatment courses that follow each other sequentially, the primary patient is denoted $p_{c,1}$ and a dummy patient is created for the secondary treatment, denoted $p_{c,2}$. The primary and secondary patients are connected in $(p_{c, 1}, p_{c, 2}) \in \mathcal{P}^{con}$. For these treatments, \eqref{Eq:seq1} states that the first day of the secondary treatment is at least one day after the last day of the primary treatment, and \eqref{Eq:seq2} states that the first day of the secondary treatment is at most three days after the last day of the primary treatment for each pair $(p_{c,1}, p_{c,2}) \in \mathcal{P}^{con}$. The full formulation is thus
\begin{alignat}{3}
	 &\text{minimize} \quad        &  & \mathmakebox[0pt][l]{1+\sum_{\pinP} \sum_{\iinK} c^i_{p} a^i_{p}}                                                                                     \label{Eq:master obj}\\
	 & \text{subject to} \quad &  &  \sum_{i \in \mathcal{K}_p} a^i_{p} = 1\,, \enskip                                                                                                                                                                      &  & \forall \pinP  \label{Eq:master a=1}\\
	 &                   &  & \sum_{\pinP} \sum_{\iinK} a^i_{p} D^i_{p,m,d,w}  +S_{m,d,w} \leq L_w, \enskip &  &  \forall \minM, \dinD, \winW         \label{Eq:master sched}      \\
	 &                 & &\sum_{i \in \mathcal{K}_ {p_{c,1}}} a^i_{p_{c,1}} \dend{p_{c,1}} + 1 \leq \sum_{i \in \mathcal{K}_ {p_{c,2}}} a^i_{p_{c,2}} \dstart{p_{c,2}},  \quad & & \forall (p_{c,1}, p_{c,2}) \in \mathcal{P}^{con} \label{Eq:seq1}\\
	 &                 & &\sum_{i \in \mathcal{K}_ {p_{c,2}}} a^i_{p_{c,2}} \dstart{p_{c,2}} \leq \sum_{i \in \mathcal{K}_ {p_{c,1}}} a^i_{p_{c,1}} \dend{p_{c,1}} + 3,  \enskip & & \forall (p_{c,1}, p_{c,2}) \in \mathcal{P}^{con}  \label{Eq:seq2}\\
	 &                   &  & a^i_{p} \in \{0,1\} \enskip & & \forall \pinP, \iinK . \label{Eq:master integer}
\end{alignat}
Relaxing the integer assumption and solving the LP yields the dual variables $\lambda_p$ associated with (\ref{Eq:master a=1}), $\gamma_{m,d,w}$ associated with (\ref{Eq:master sched}), $\eta_c$ associated with \eqref{Eq:seq1} and $\xi_c$ associated with \eqref{Eq:seq2}, where each $c$ is the index of the $(p_{c,1},p_ {c,2})$-pair.

\subsection{Subproblems: Schedule Generation}
\label{Sec:subproblems}
One subproblem is formed for each patient $\pinP$, with the aim to generate a new feasible schedule to add to $\mathcal{K}_p'$, i.e., as a column to the RMP. The main variables for the subproblem model are presented in Table~\ref{tab:variables_IP}. As the subproblems are complex, the constraints and objectives are described individually.

\begin{table}[ht]
\caption{Main variables in the subproblem (schedule generation)}
\label{tab:variables_IP}
\begin{small}
\begin{tabular}{P{0.15\linewidth} P{0.798\linewidth}}
\hline\noalign{\smallskip}
	$\q{d}{f} \in \{0,1\}$                 & $1$ if  fraction $\finF_p$ is scheduled on weekday $\dinD$ on machine $\minM$ in schedule $\iinK$ for $\pinP$, $0$ otherwise\\
		$x^i_{p,m,d,w} \in \{0,1\}$             & $1$ if patient $p$ in schedule $\iinK$ is scheduled in window $\winW$ on machine $\minM$ on weekday $\dinD$, $0$ otherwise \\
		$t^i_{p,m,d,w} \in \{0,1\}$             & $1$ if patient $p$ in schedule $\iinK$ starts treatment in window $\winW$ on machine $\minM$ on weekday $\dinD$, $0$ otherwise \\
    $\nu^i_{p,g} \in \{0,\dots, 5\}$                & The number of fractions scheduled in week $\ginG$ for patient $\pinP$ in schedule $\iinK$\\
    $\tau^i_p \in \{1,2,\dots \}$                  &The total number of weeks that each patient is scheduled  \\
     $\rho^i_{p,g} \in \{0,1\}$             &0 for week $g$ if the minimum fraction requirement is not met, or if $g$ is the first or last week of treatment \\
 \noalign{\smallskip}\hline
\end{tabular}
\end{small}
\end{table}

\subsubsection{Subproblem Constraints}
The subproblem constraints should ensure that all medical and technical constraints are fulfilled in schedule $\iinK$ for each patient $\pinP$.
Constraint \eqref{Eq:sum(q) over d,m =1} forces fraction $f$ to be scheduled exactly one time for each patient. Constraint \eqref{Eq:q_1 = t} states that the first fraction for patient $p$ is scheduled on machine $m$ on day $d$, in any window, whereas constraint \eqref{Eq:t <= x} also gives the correct time window $w$ for the first fraction. Furthermore, constraint \eqref{Eq:sum(x) = sum(q)} states that each patient is scheduled in exactly one time window for each fraction.
\begin{alignat}{3}
& &&\sum_{\minM} \sum_{\dinD} \q{d}{f} = 1, \enskip           &  & \forall  \finF_{p}                                       \label{Eq:sum(q) over d,m =1} \\
& && q^i_{p,m,d,1} = \sum_{\winW} t^i_{p,m,d,w} , \enskip               &  & \forall \minM, \dinD                                   \label{Eq:q_1 = t}           \\ 
& &&  t^i_{p,m,d,w} \leq x^i_{p,m,d,w} , \enskip                                   &  & \forall  \minM, \dinD  , \winW                   \label{Eq:t <= x}   \\
& &&  \sum_{\winW} x^i_{p,m,d,w} = \sum_{\finF_p} \q{d}{f},  \quad       &  & \forall \minM, \dinD                \label{Eq:sum(x) = sum(q)}   
 \end{alignat}
The earliest day to start treatment is $\dMin$ and a treatment can only start on an allowed start day given by $\Ahp$. The fractions can only be scheduled on a machine for the patient protocol given by $\mathcal{M}_{h_p}$, but not if the machine is unavailable. In total, this is captured in constraints \eqref{Eq:q = 0 daylim} and \eqref{Eq:q = 0 M_p}.
\begin{alignat}{3}
& && q^i_{p,m,d,1}  = 0,  \quad    &  & \forall \minM, \dinD \text{ if } d<\dMin \text{ or } d \notin \Ahp \label{Eq:q = 0 daylim}         \\
 & && \q{d}{f}  = 0,  \quad    &  &\forall \minM, \dinD \text{ if } m \notin \mathcal{M}_{h_p} \text{ or } d \in \DmU, \finF_p      \label{Eq:q = 0 M_p}        
\end{alignat}

Constraint \eqref{Eq:colgen sub sched} ensures that the treatment fits within each time window  $w$ on machine $m$ on day $d$. Since the duration of the first fraction is different from the rest, the first term will evaluate to zero in the first fraction. The duration plus the already occupied time slots $S_{m,d,w}$ in that window should be less than or equal the window length $L_w$. 
 \begin{align}
 \Big((x^i_{p,m,d,w} - t^i_{p,m,d,w}) \durprest +  t^i_{p,m,d,w}\durpfirst \Big) +  S_{m,d,w} \leq L_w,  \quad                                                                                                              &  \forall \minM, \dinD, \winW                                      \label{Eq:colgen sub sched}   
 \end{align}
 

\mypar{Planned machine unavailability}
Because of planned machine unavailability due to maintenance or holidays, it is not possible to always schedule fractions on consecutive weekdays the way it was done in \cite{Frimodig2022}. Instead, fractions are scheduled on consecutive days on weeks with no unavailability in machine group by constraint \eqref{Eq:q(d,f)=q(d+1,f+1) on consecutive day periods}. Furthermore, two fractions are scheduled on consecutive days when two days in a row are available for all machines in a machine group, for patients where the minimum fractions per week is less than five (otherwise, two fractions on the same day may be necessary if there is an unavailable day is in the same week) in \eqref{Eq:q(d,f)=q(d+1,f+1) on consecutive days}. All fractions must be in the same machine group for all patients, which is enforced by constraint \eqref{Eq:Machine group}. Finally, since the duration is different for the first fraction, and because the day of the last fraction is needed in the objectives, constraint \eqref{Eq:f_order} enforces an ordering of the fractions. 
\begin{alignat}{3}
&  &&\begin{multlined} \sum_{\minbM} \q{d}{f} = \sum_{\minbM} \q{d+1}{f+1}, \end{multlined} \enskip && \begin{multlined} \forall  \binBm, 
 \dinD_{g}, \finF_p, \ginG \text{ where } g \notin \Gubm\end{multlined}  \label{Eq:q(d,f)=q(d+1,f+1) on consecutive day periods} \\
& &&\begin{multlined} \sum_{\minbM} \q{d}{f} = \sum_{\minbM} \q{d+1}{f+1}, \\ \end{multlined} \enskip && \begin{multlined} \text{if } \fminHp < 5, \forall \binBm, \finF_p, \dinD \\  \text{where } d \notin \DmU \text{ and } d+1 \notin \DmU \text{ for all } m \in b_\mathcal{M}\end{multlined} \label{Eq:q(d,f)=q(d+1,f+1) on consecutive days} \\
& &&\sum_{\minbM} \sum_{\dinD} \q{d}{f} = \sum_{\minbM} \sum_{\dinD}  \q{d+1}{f+1}, \quad && \forall \binBm, \finF_p
\label{Eq:Machine group} \\
& && \sum_{\dinD} d \sum_{\minM} \q{d}{f} \leq \sum_{\dinD} d \sum_{\minM} \q{d}{f+1}, \quad && \forall f=\{1,\dots, F_p-1\}
\label{Eq:f_order}
\end{alignat}
If all machines in a machine group is available a week, there should be at most one fraction per day, enforced by \eqref{Eq:Max1frac}. If a patient requires less than five fractions per week, constraint \eqref{Eq:Max1frac_2} states that there should also be at most one fraction per day. If there is unavailability in a week, there should instead be at most two fractions on one day, enforced by \eqref{Eq:Max2frac}. Constraint \eqref{Eq:pairwise3} states that at most one day per week has two fractions scheduled on the same day. There should be at most two gap days between fractions, stated in constraint \eqref{Eq:pause_days}. Finally, it is not allowed to have two fractions on the first day of treatment, which is enforced by constraint \eqref{Eq:first_day}.
\begin{alignat}{3}
& &&\sum_{\minbM} \sum_{\finF_p} \q{d}{f} \leq 1, \quad &&\forall \binBm, \dinD_g, \ginG \text{ where } g \notin \Gubm \label{Eq:Max1frac}\\
& &&\sum_{\minM} \sum_{\finF_p} \q{d}{f} \leq 1, \quad && \text{if } \fminHp < 5, \forall \dinD \label{Eq:Max1frac_2}\\
& &&\sum_{\minbM} \sum_{\finF_p} \q{d}{f} \leq 2, \quad &&\text{if } \fminHp = 5, \forall \binBm, \dinD_g, g \in \Gubm \label{Eq:Max2frac} \\
& &&\sum_{\minM} \sum_{\finF_p} (\q{d_1}{f}+ \q{d_2}{f}) \leq 3, \quad &&\forall  \ginG, d_1, d_2 \in \mathcal{D}_g \text{ with } d_1 < d_2  \label{Eq:pairwise3}\\
& &&\sum_{\dinD} d  \sum_{\minM} (\q{d}{f+1}-\q{d}{f}) \leq 3,  \quad && \forall f=\{1, \dots, F_p-1 \}\label{Eq:pause_days} \\
& && \sum_{\dinD} d \sum_{\minM}  ( \q{d}{2} - \q{d}{1} )\geq  1 \quad && \text{if } F_p > 1 \label{Eq:first_day}
\end{alignat}
If two fractions are scheduled on the same day, they need to be scheduled with some time apart (typically 6 hours). This is enforced by constraint \eqref{Eq:dd_1} (note that $\mathcal{W}=\{1,2,3,4\}$). If there are no treatments scheduled on day $d$, \eqref{Eq:dd_1}  simply states that $0\leq2$. If there is one fraction scheduled, the fraction can be scheduled in any window that day. If there are two fractions scheduled, the the fractions cannot be scheduled in window two or three since the left-hand-side will then be greater than two, forcing them to be scheduled in windows one and four, and thereby being far enough apart. Furthermore, constraint \eqref{Eq:dd_one_per_window} forces maximum one fraction to be scheduled in each time window.
\begin{alignat}{3}
& &&2\sum_{\minM} \sum_{w=\{2,3\}} x^i_{p,m,d,w} + \sum_{\minM} \sum_{w=\{1,4\}} x^i_{p,m,d,w} \leq 2  , \quad &&\forall \dinD \label{Eq:dd_1}\\
& &&\sum_{\minM} x^i_{p,m,d,w} \leq 1\quad     &&\forall \winW, \dinD \label{Eq:dd_one_per_window}
\end{alignat}

\mypar{Minimum Fractions per Week}
The minimum fractions per week requirement does not apply in the weeks before treatment starts or after it ends. It is approved to have less than the minimum fractions in the first or last treatment week, since a treatment could start or end in the middle of the week. The minimum number of fractions per week must be fulfilled in all intermediate weeks between the first week and the last. The number of fractions that are scheduled in week $g$ are counted in the variable $\nu^i_{p,g}$ in equation \eqref{Eq:nu}. The total number of weeks that each patient is scheduled is computed as the variable $\tau^i_p$ by subtracting the start week from the end week in \eqref{Eq:tau}. This gives the number of intermediate weeks, that have minimum fraction requirements. For patients with $\fminHp < 5$, there can be intermediate weeks if there are six (or more) fractions in total. Therefore, this constraint applies for all patients where $F_p > 5$. 
\begin{alignat}{3}
& && \nu^i_{p,g} = \sum_{\minM} \sum_{\dinD_g} \sum_{\finF_p} \q{d}{f}, \quad && \forall \ginG \label{Eq:nu} \\
& && \tau^i_p = \sum_ {\ginG} g \sum_{\dinD_g} \sum_{\minM} ( \q{d}{F_p}  - \q{d}{1} ) - 1, \quad && \text{if } F_p > 5 \label{Eq:tau}
\end{alignat}
If the minimum fraction requirement is not met, i.e., the number of scheduled fractions $\nu^i_{p,g} < \fminHp$, then \eqref{Eq:rho} forces $\rho^i_{p,g} = 0$. Also, $\rho^i_{p,g} = 0$ if the week is the start week or the end week by constraints \eqref{Eq:rho_startweek} and \eqref{Eq:rho_endweek}. 
\begin{alignat}{3}
& &&\nu^i_{p,g} \geq \fminHp \rho^i_{p,g}, \quad &&\forall \ginG \label{Eq:rho} \\
& &&\rho^i_{p,g} + \sum_{\minM} \sum_{\dinD_g}  \q{d}{1}  \leq 1, \quad  &&  \text{if } F_p > 5,\forall \ginG \label{Eq:rho_startweek} \\
& &&\rho^i_{p,g} + \sum_{\minM} \sum_{\dinD_g} \q{d}{F_p}  \leq 1, \quad  && \text{if } F_p > 5,\forall \ginG  \label{Eq:rho_endweek}
\end{alignat}
All intermediate weeks must have minimum fraction requirement scheduled. This is enforced by summing all binary variables $\rho^i_{p,g}$ to be the number of intermediate weeks in constraint \eqref{Eq:rho_tau}. 
\begin{align}
\sum_{\ginG} \rho^i_{p,g} = \tau^i_p, \quad  \text{if } F_p > 5 \label{Eq:rho_tau}
\end{align}

\mypar{Non-Conventional Treatments}
The most common type on non-conventional treatment is where patients are treated with a pause days between each fraction. For these patients, constraints \eqref{Eq:q(d,f)=q(d+1,f+1) on consecutive day periods}, \eqref{Eq:q(d,f)=q(d+1,f+1) on consecutive days} (fractions on consecutive days) do not apply. 
Since the time horizon is assumed to have only weekdays, but also the weekend is considered a pause between fractions, the set $\mathcal{D}^{\textit{Fri}}$ is created to include all Fridays in the planning horizon. Constraint \eqref{Eq:non-conv, 2 days apart} states that fraction $f$ and $f+1$ should be at least two days apart, unless it is a Friday and then \eqref{Eq:non-conv, 1 days apart} states that the fractions should instead be at least one weekday apart. Constraint \eqref{Eq:pause_days} already states that the fractions should be at most three days apart.
\begin{alignat}{3}
& && \sum_{\dinD \setminus \mathcal{D}^{\textit{Fri}}} d \sum_{\minM} ( \q{d}{f+1} - \q{d}{f} ) \geq 2 \quad && \forall f=\{1,\dots,F_p-1\} \label{Eq:non-conv, 2 days apart} \\
& && \sum_{d \in \mathcal{D}^{\textit{Fri}}} d \sum_{\minM} ( \q{d}{f+1} - \q{d}{f} ) \geq 1 \quad &&  \forall f=\{1,\dots,F_p-1\} \label{Eq:non-conv, 1 days apart} 
\end{alignat}
Another type of non-conventional treatment is for a group of patients where $F_p=5$ and $\fminHp=5$, meaning that all fractions must be scheduled in the same week Monday to Friday. For these patients, $\mathcal{A}_p$ is the set of Mondays in the planning horizon, and the ordering is enforced by constraint \eqref{Eq:5perweek}.
\begin{align}
\sum_{\minbM} \q{d}{f} = \sum_{\minbM} \q{d+1}{f+1}, \quad \forall \binBm, d=\{1,\dots,D-1\}, f=\{1,\dots,F_p-1\} \label{Eq:5perweek}
\end{align}
A third type of non-conventional treatment at Iridium Netwerk is total body irradiation (TBI) treatments, that should be delivered twice daily on Monday to Wednesday. For these patients, $\mathcal{A}_p$ is the set of Mondays in the planning horizon and $\fminHp=6$. Constraints \eqref{Eq:Max1frac} and \eqref{Eq:Max1frac_2} are removed to allow multiple fractions per day. Constraint \eqref{Eq:dd_1} already forces the fractions to be scheduled in the first and last time window, but we now allow this for multiple days per week by removing constraint \eqref{Eq:pairwise3}. Constraints \eqref{Eq:q(d,f)=q(d+1,f+1) on consecutive day periods}, \eqref{Eq:q(d,f)=q(d+1,f+1) on consecutive days}  and \eqref{Eq:first_day} are also removed for these patients to never force fractions to be on consecutive days. Furthermore, \eqref{Eq:non-conv, TBI} states that there should be at most two days between first and last fraction. 
\begin{align}
\sum_{\dinD}d \sum_{\minM} (\q{d}{F_p} -\q{d}{1}) = 2 \label{Eq:non-conv, TBI}
\end{align}

\subsubsection{Subproblem Objective Function} 
The subproblem objective function has two main components; the cost for the schedule, and the cost related to the dual variables from the master problem. The variables used to formulate the objective function are presented in Table~\ref{tab:variables_sub_obj}.
\begin{table}[ht]
\caption{Cost function variables in the subproblem (schedule generation)}
\label{tab:variables_sub_obj}
\begin{small}
\begin{tabular}{P{0.16\linewidth} P{0.78\linewidth}}
\hline\noalign{\smallskip}
	$z^i_{p,d} \in \{0,1\}$                       & $1$ if patient $p$ in schedule $i$ has switched windows from day $d$ to $d+1$, $0$ otherwise\\
	$u^i_{p} \in \{0,1, \dots\}$                   &The number of violation of the time window preference for patient $\pinP$ in schedule $i$ \\
	$s^i_{p,f} \in \{0,1\}$                   & $1$ if patient $p$ in schedule $i$  switches to a partially beam-matched machine from fraction $f$ to $f+1$, $0$ otherwise, for $f=\{1,\dots,F_p-1\}$   \\
	 $o^i_p  \in \{0,1,\dots\}$         & Total excess time in days for patient $p$ in schedule $i$, meaning the number of days the treatment has been paused\\
 \noalign{\smallskip}\hline
\end{tabular}
\end{small}
\end{table}
\mypar{Schedule cost} The objectives presented in Section \ref{Sec: Problem Formulation} are formulated as a cumulative cost function, where the different objectives are combined with weights $\alpha_1,\alpha_2,\alpha_3, \alpha_4, \alpha_5, \alpha_6, \alpha_7$ in \eqref{Eq:colgen sub cost} to the total cost $c^i_{p}$ of the schedule $\iinK$.
\begin{align}
\begin{split}
c^i_{p} =\enskip &\alpha_1 ( c_{h_p} \sum_{\minM} \sum_{d=\dMin}^{D} q^i_{p,m,d,1}(d-\dMin) )+  \alpha_2( \sum_{\dinD} z^i_{p,d} ) + \\
& \alpha_3 (u^i_{p}) +
\alpha_4 ( \sum_{\finF_p} \sum_{\minM} \sum_{\dinD} \q{d}{f} \mathbb{1}_{(m \notin \mathcal{M}^{\mathit{pref}}_p)}) + 
 \alpha_5 (\sum_{f=1}^{F_p-1} s^i_{p,f}) + \\
 & \alpha_6 (\sum_{\finF_p} \sum_{\minM} \sum_{\dinD} \q{d}{f} \mathbb{1}_{(m \notin  \SMpref)} ) + 
\alpha_7 o^i_{p}
\label{Eq:colgen sub cost}
\end{split}
\end{align}
The most important objective is to minimize a weighted sum of the waiting times. In the first term in \eqref{Eq:colgen sub cost}, the number of waiting days after $\dMin$ are linearly penalized with weight $c_{h_p}$ corresponding to the priority group of protocol $h$ for patient $p$.
The second cost is the deviations in treatment time for each patient, which is computed as the number of time window switches. 
The variable $z_{p,d}$ is defined according to constraints \eqref{Eq:MIP_z_1} and \eqref{Eq:MIP_z_2} to compute the time window switches between two days, and used in the second term in the cost function  \eqref{Eq:colgen sub cost}. 
\begin{align}
z^i_{p,d}& \geq \sum_{\minM} (x^i_{p,m,d,w}-x^i_{p,m,d+1,w}), \enskip              &  & \forall  d=\{1,\dots,D-1\}, \winW                                     \label{Eq:MIP_z_1}\\
z^i_{p,d} &\geq \sum_{\minM} (x^i_{p,m,d+1,w}-x^i_{p,m,d,w}), \enskip              &  & \forall d=\{1,\dots,D-1\}, \winW                                       \label{Eq:MIP_z_2} 
\end{align}
Some patients have a preference on treatment time of the day. The variable $u^i_{p}$ is defined in \eqref{Eq:MIP_u1} and \eqref{Eq:MIP_u2} to be the violation of the time window preference for each patient, and is used in the third term in \eqref{Eq:colgen sub cost}. 
\begin{align}
 u^i_{p} &= \sum_{\minM}  \sum_{\dinD}  \sum_{\winW} x^i_{p,m,d,w}\mathbb{1}_{(w \neq \wPref)}, \enskip      &  & \text{if }\pinPpref \label{Eq:MIP_u1} \\
u^i_{p} &= 0, \enskip   &  & \text{if }p \notin \Ppref                             \label{Eq:MIP_u2} 
\end{align}
The number of fractions scheduled on a non-preferred machine stated by the treatment protocol are summed in the fourth term in \eqref{Eq:colgen sub cost}.
Moreover, there is a cost for the number of switches between machines that are not completely beam-matched. If fraction $f$ is scheduled on a machine in a group of completely beam-matched machines, but $f+1$ is not, then it must be scheduled on a partially matched machine. The variable $s^i_{p,f}$ is one if there is a switch to a partially matched machine, enforced by constraint \eqref{Eq:IP_machineswitch}. All machine switches to partially matched machines are summed in the fifth term in \eqref{Eq:colgen sub cost}.
\begin{align}
s^i_{p,f}&\geq \sum_{\dinD} \sum_{m \in c_\mathcal{M}} (q^i_{p,m, d, f} - q^i_{p,m, d, f+1}), \quad  \forall c_\mathcal{M} \in \mathcal{C_M}, f = \{1,\dots,F_p-1\} \label{Eq:IP_machineswitch}
\end{align}
Most patients have preferences on what hospital to be treated on based on where they live; $\SMpref$ is the list of machines at the preferred site of patient $\pinP$. In the sixth term  in \eqref{Eq:colgen sub cost}, the fractions scheduled on another site than the preferred is summed. 
Finally, there is an objective to keep the overall treatment time as short as possible. The overall treatment time is computed as the number of weekdays from the first to the last fraction. The excess time has a lower bound of zero, and is otherwise the number of days extra needed to complete the fractions, computed by \eqref{Eq:tot_overdue_time} and added to the cost in the seventh term in \eqref{Eq:colgen sub cost}.
\begin{align}
o^i_{p} &\geq \sum_{\dinD} d \sum_{\minM}  (q^i_{p,m,d,F_p} - q^i_{p,m,d,1}) - F_p  \label{Eq:tot_overdue_time}
\end{align}
The values of the weights $\alpha_1, \dots, \alpha_7$ should reflect the importance of the objectives in relation to each other. In \cite{Frimodig2022}, it was shown that the balancing of the different objectives is clearly related to the weights. In this paper, the weights are fixed to have values $\alpha_1=100,\alpha_2=1,\alpha_3=1, \alpha_4=10, \alpha_5=10, \alpha_6=50, \alpha_7=300$. Since waiting time objective is also weighted by $c_{h_p}$, this objective has the highest weight for priority A patients.

\mypar{Objective function} Relaxing the integer assumption in the master problem and solving the LP yields the dual variables $\lambda_p$ associated with (\ref{Eq:master a=1}), $\gamma_{m,d,w}$ associated with (\ref{Eq:master sched}), $\eta_c$ associated with \eqref{Eq:seq1} and $\xi_c$ associated with \eqref{Eq:seq2}. The master problem constraints have parameters $D^i_{p,m,d,w}$, corresponding to the duration of the treatment on machine $m$, day $d$, time window $w$ in schedule $\iinK$, and $\dstart{p}$ and $\dend{p}$ stating the first and last day of the treatment in schedule $\iinK$. These parameters are computed from the subproblem variables in \eqref{Eq:Dur_master}, \eqref{Eq:dstart} and \eqref{Eq:dend}.
\begin{align}
D^i_{p,m,d,w} &= \Big((x_{p,m,d,w}^i - t_{p,m,d,w}^i) \durprest  + t_{p,m,d,w}^i \durpfirst\Big) \quad \forall \minM, \dinD, \winW \label{Eq:Dur_master} \\
\dstart{p} &= \sum_{\dinD} d \sum_{\minM} \q{d}{1} \label{Eq:dstart}\\
\dend{p} &= \sum_{\dinD} d \sum_{\minM} \q{d}{F_p} \label{Eq:dend}
\end{align}
For the consecutive treatments, the subproblem objective function has an additional term: for primary treatments, the term is defined in \eqref{Eq:obj_seq_p1}, and for secondary treatments, the term is defined in \eqref{Eq:obj_seq_p2}, both based on the dual variables from \eqref{Eq:seq1} and \eqref{Eq:seq2} with $\dstart{p}$ and $\dend{p}$ defined in \eqref{Eq:dstart}, \eqref{Eq:dend},  where $c$ is the index of the $(p_{c,1},p_ {c,2}) \in \mathcal{P}^{con} $-pair. 
\begin{align}
 -(\eta_c-\xi_c) \sum_{\dinD} d   \sum_{\minM} \q{d}{F_p} \quad  &\text{ if } p = p_{c,1} \text{ for } (p_{c,1},p_{c,2})  \in \mathcal{P}^{con}   \label{Eq:obj_seq_p1} \\
 -(\xi_c-\eta_c) \sum_{\dinD} d   \sum_{\minM} \q{d}{1}  \quad& \text{ if } p = p_{c,2}  \text{ for } (p_{c,1},p_{c,2})  \in \mathcal{P}^{con}  \label{Eq:obj_seq_p2}
\end{align}

The subproblem objective function \eqref{Eq:IP objective} is the cost of the schedule defined by \eqref{Eq:colgen sub cost}, minus the master dual variables multiplied by the coefficients given from their respective constraints in the master problem. If $p \in\mathcal{P}^{con}$, the term from \eqref{Eq:obj_seq_p1} or \eqref{Eq:obj_seq_p2} is also added.
\begin{align}
	  \text{minimize} \quad   
	 c^i_{p} -  \lambda_p - \sum_{\minM} \sum_{\dinD} \sum_{\winW} \gamma_{m,d,w} \Big((x_{p,m,d,w}^i - t_{p,m,d,w}^i) \durprest
+  t_{p,m,d,w}^i \durpfirst \Big)  \label{Eq:IP objective} 
\end{align}

\subsection{Heuristic to Create Initial Schedules}
\label{Sec:CG-heuristic}
When solving a CG formulation, it is often beneficial to have a heuristically generated set of initial columns \cite{Vanderbeck2005}, which is true also for the CG algorithm as described in Figure~\ref{Fig:ColGenFlow}. For the algorithm to perform well, it is essential that the set of initial schedules are of good \emph{quality}, and that there is \emph{diversity} in the initial schedules for each patient (i.e., not the same column generated over again). Therefore, there must be randomization included in the heuristic. The number of initial schedules is set to 50, because for a larger number it becomes difficult to create new schedules that are not duplicates of the ones already created without sacrificing too much in quality, which in general increases solution times. 

The overarching idea with the schedule construction algorithm is that it will create 50 \emph{complete} schedules with all patients, that are varying in quality but that all fulfill the capacity constraints. In other words, for $i=1,\dots,50$ the capacity constraint in the master problem \eqref{Eq:master sched} can be fulfilled by setting $a^i_p=1$ for all $\pinP$ for one $i$ at a time. This is achieved in the algorithm by looping through all patients (in partly randomized order), and assigning them to machines, days and time windows according to some specific (partly randomized) order, while making sure that the maximum capacity of the resources is not exceeded, and repeating this from scratch for each $i$. The time horizon is chosen long enough so that it will always be possible to find feasible schedules on some machine. The schedule construction heuristic is presented in pseudocode in Appendix \ref{Appendix:heuristic}. 

\section{Uncertainty in Future Arrivals}
\label{Sec:FutureArrivals}
\begin{figure}[b]
\center
\includegraphics[scale=0.9]{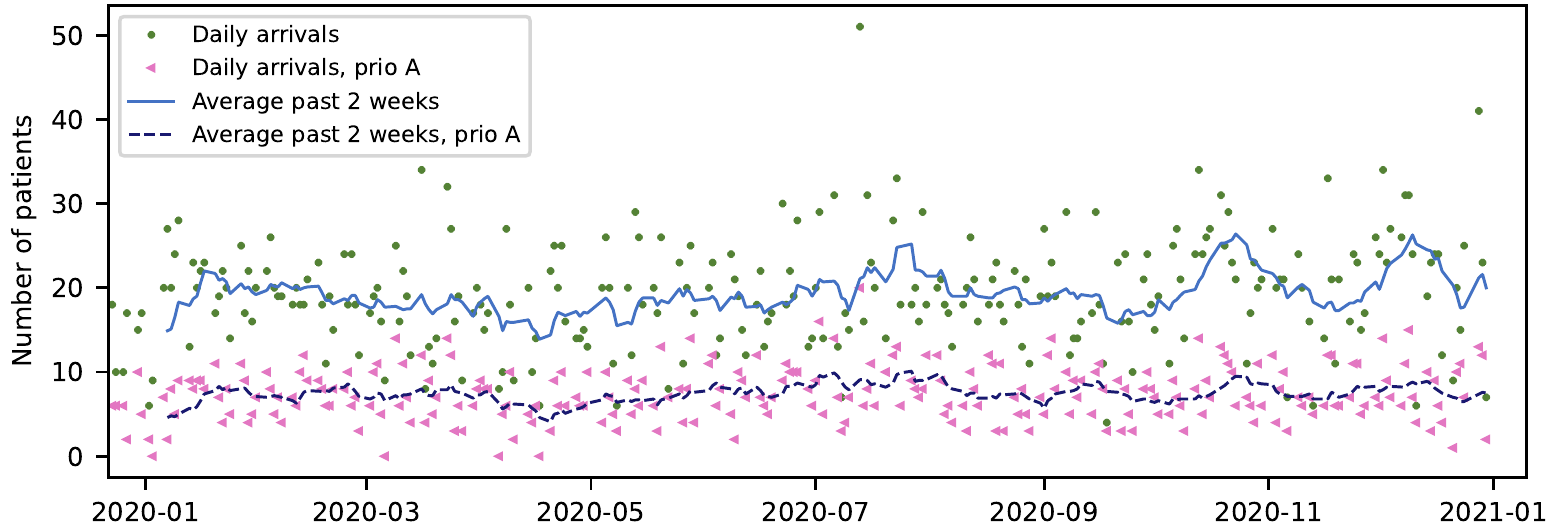}
\caption{Daily arrivals in 2020 at Iridium Netwerk, excluding weekends and holidays}
\label{Fig:Arrivals}
\end{figure}

A major challenge in RT patient scheduling is the stochastic patient inflow, and that patients of high priority should start treatment as soon as possible after arrival. Figure~\ref{Fig:Arrivals} shows the daily patient arrivals at Iridium in 2020. The two-week rolling average varies between 13.9 and 26.4 patients per day, and between 4.1 and 10.1 priority A-patients per day. 
Most clinics, including Iridium Netwerk, reserve a proportion of machine capacity for high priority patients each day. This solution can result in poor quality schedules; urgent patients may have to wait for treatment if not enough capacity was reserved, or the treatments for low-priority patients may be delayed as a result of poor machine utilization. 

To account for uncertainties in future arrivals, three different methods to reserve time for high priority patients are investigated. These include to not save time for future arrivals at all (\emph{no reservation}), to reserve time on the machines based on average utilization rates for high priority patients, which should mimic the current practices (\emph{static time reservation}), and to include expected future patients arrivals in the problem as placeholder (dummy) patients, thereby allowing for trade-offs with the actual patients (\emph{dynamic time reservation}). 
Assuming that we could see into the future, it is possible to do an \emph{offline} solution by using the actual future arrivals for the coming weeks, and including them each day when creating the schedule. This is of course not possible in reality, butis used to get a best case schedule.

\begin{figure}
\center
\includegraphics[scale=0.9]{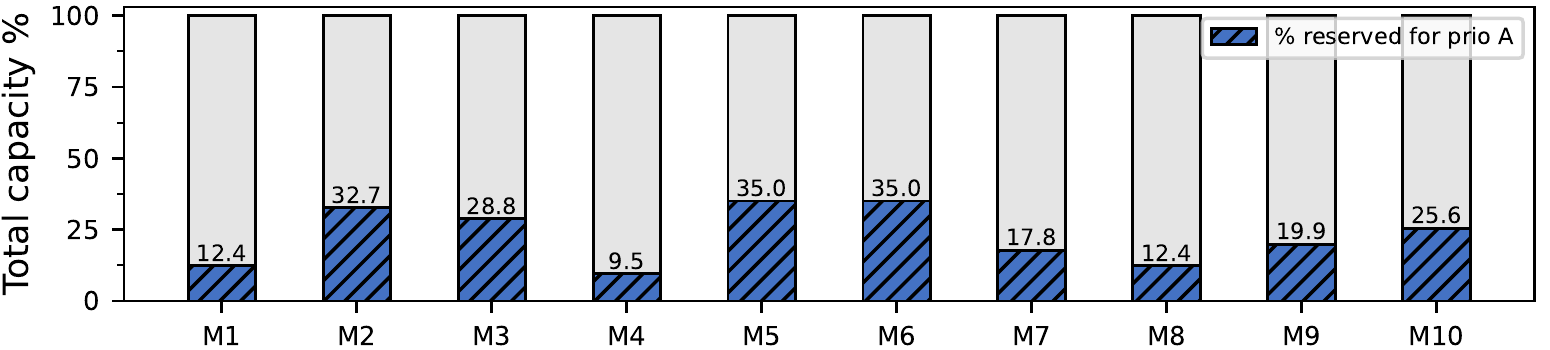}
\caption{Percentage capacity reserved for priority A patients for each machine for \emph{static time reservation} method}
\label{Fig:static_reserved}
\end{figure}

The \emph{static time reservation} method uses the historical patient arrivals, and distributes them over the preferred machines for each protocol based on the average number of fractions and average session time for each protocol according to Figure~\ref{Fig:static_reserved}. The time is reserved for priority A patients by blocking it for priority B and C patients; the input schedule $S_{m,d,w}$ is adjusted to include the blocked time when schedules for priority B and C patients are generated, both in Algorithm~\ref{Alg:ColGenGreedy} and constraint \eqref{Eq:colgen sub sched}. However, the right-hand-side in \eqref{Eq:master sched} remains intact, since the constraint deals with patients of all priority groups. This is not a problem since no schedules that uses the blocked time can be created.

The \emph{dynamic time reservation} method adds the expected future priority A patients as \emph{placeholder} (dummy) patients and includes them in the scheduling. Based on historical arrival rates, 36 priority A patients are expected to arrive each week on average. In the dynamic reservation method, these are distributed over the protocols according to the protocol probabilities, where the protocols with similar characteristics are grouped to simplify the problem.  An overview of the placeholder patients for each week can be seen in Table~\ref{Tab:placeholder_future}. In each daily scheduling problem, the placeholder patients are added for each week in the planning horizon where the actual patients are expected to start their treatments, which is usually around 8 weeks. The cost function $c^i_{p}$ in \eqref{Eq:colgen sub cost} for the future dummy patients does not include costs for time window switches, time window preferences, or machine switches, but it does include waiting time, preferred machines, overall treatment time and site preferences. The site preferences are generated for each dummy future patient based on the historical distribution between the hospitals. To include the expected future patients as placeholder patients in the scheduling will reserve capacity on machines for future patients, however, the reservation is dynamic since it allows for trade-offs with the actual patients; if a low-priority patient has already waited a very long time for treatment, there are cases where it is beneficial for the overall schedule to let that patient start before a patient of higher priority. 
\begin{table}
\caption{Placeholder priority A patients added to problem in the dynamic time reservation method}
\label{Tab:placeholder_future}
\setlength\extrarowheight{2pt}
\begin{tabular}{P{1.3cm}P{1.5cm}P{2.1cm}P{1.95cm}P{3.4cm}P{1.3cm}P{2.0cm}}
\hline
Protocol & \# patients/ week & Minimum fractions/week & $\{ \durpfirst,  \durprest\}$ & Machines preferred & Number fractions & Minimum days from CT\\
\hline
Urgent 1 & 19 & 1 & $\{24, 24\}$ & M1, M2, M3, M4, M5, M6, M7, M8 & 3 &  0 \\
VMAT 1 & 6 & 5 & $\{24, 12\}$ & M2, M3, M5, M6, M10 & 28 & 11 \\
VMAT 2 & 5 & 4 &$\{24, 12\}$ & M1, M2, M3, M5, M6, M7, M8, M10 & 23 & 10 \\
STX 1     & 3 & 3 & $\{40, 40\}$ & M9 & 6 & 5 \\
Electron & 1 & 3 & $\{24, 12\}$ &M1, M4, M5, M6, M8 & 12 & 5 \\
VMAT 3 & 1 & 3 & $\{24, 12\}$ & M1, M4, M5, M6, M8 & 10 & 4 \\
Urgent 2 & 1 & 1 & $\{24,24\}$ & M9 & 1 & 0 \\ 
\hline
\end{tabular}
\end{table}

\section{Results}
\label{Sec:results}

\begin{figure}[ht]
\center
\includegraphics[scale=0.9]{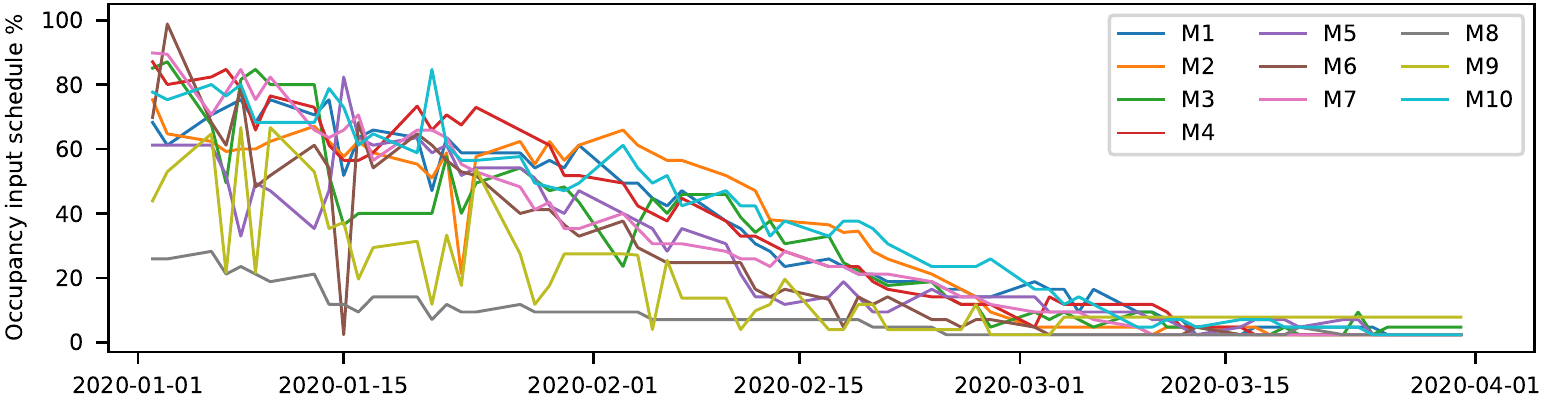}
\caption{Occupancy in input schedule resulting from patients that arrived in 2019}
\label{Fig:Input_sched_occupancy}
\end{figure}

To test the model, data from Iridium Netwerk from 2020 is used. In 2020, they operated ten linacs. One of them, M9, is a specialized linac and its scheduling can be done separately from the other nine linacs. In 2020, there were 87 days with planned machine unavailability, not counting the weekends. This means that around 34\% of all weekdays in the year had some type of machine unavailability planned due to holidays, maintenance or other quality assurance activities.  
All patients that arrived to Iridium Netwerk in 2019, but had treatments scheduled in 2020, are seen as fixed in the input schedule for the 2020 scheduling. The occupancy in the 2020 schedule resulting from 2019 patients can be seen in Figure~\ref{Fig:Input_sched_occupancy}. On the first day of 2020, this occupancy will make up the input schedule $S_{m,d,w}$ in the CG model. 

The daily batch scheduling at Iridium Netwerk is simulated starting from the 1\textsuperscript{st} of January 2020 by using the \emph{actual} patient arrivals at Iridium Netwerk (Figure~\ref{Fig:Arrivals}) for each day of 2020. A problem instance is made up of the \emph{input schedule} containing the treatments that are fixed to the schedule due to previous scheduling decisions, together with the list of unscheduled patients \emph{from previous days}, and the \emph{current day's arrivals}. The CG algorithm is run to perform the daily batch scheduling. In the resulting schedule, the patients that start treatment within the notification period are fixed to the input schedule for the next day, while the scheduling decisions are postponed to the next day for the remaining patients.

In total, there are 254 problem instances corresponding to the workdays in 2020. The planning horizon is three months. The number of patients to be scheduled each day varies depending on what uncertainty method from Section~\ref{Sec:FutureArrivals} that is used, but for the static time blocking method the number of patients will vary between 27 and 122 with an average of 72.3 patients per daily batch scheduling instance.
All data used in this paper is publicly available\footnote{Access through this link: \url{https://osf.io/j2bxp/?view_only=e1402382b67f4ad0a4b8a3f4ed28088a}}.
The numerical experiments are run on a Windows 10 computer, with an  Intel\textsuperscript{\textregistered} Core{\texttrademark} i9-7940X X-series processor and 64~GB of RAM. The CG model is created in Python 3.8 and solved using IBM ILOG CPLEX 20.1 in the Python API.

\subsection{Planned machine unavailability}

\begin{figure}[ht]
 \center
 \includegraphics[scale=0.9]{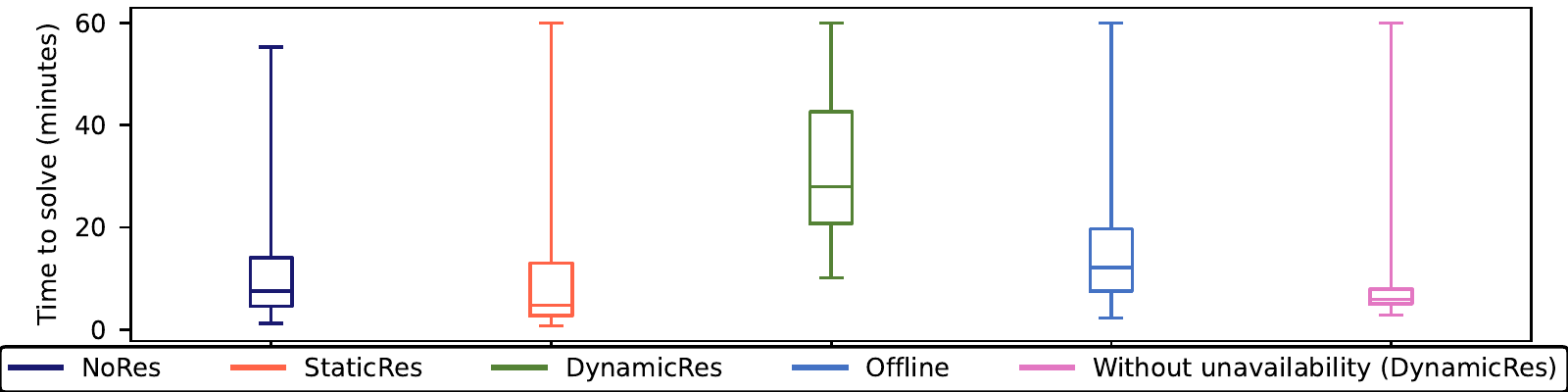}
 \caption{Time to solve the instances for the different methods to handle unavailability. "Res" here means reservation. The dynamic time reservation is shown both with and without unavailability in model. Timeout was set to 1 hour}
  \label{Fig:with_without_unavailability}
 \end{figure}

To the best of our knowledge, the RTSP model presented in this paper is the first to consider planned unavailability on machines. Without the planned unavailability, constraints \eqref{Eq:q(d,f)=q(d+1,f+1) on consecutive day periods} to \eqref{Eq:rho_tau} can be replaced by a single constraint forcing fractions to be on consecutive days, as done in \cite{Frimodig2022}. Furthermore, since gaps in the schedules induced by the scheduling procedure are only allowed when there is machine unavailability, the minimum fractions per week requirements are unnecessary.  Figure~\ref{Fig:with_without_unavailability} shows boxplots for the solution times for the different methods to handle future uncertainty. The solution times  for the the dynamic time reservation model are shown with and without planned unavailability. The time limit was set to one hour. 
One can see in Figure~\ref{Fig:with_without_unavailability} that including the planned machine unavailability in the CG model increases the solution times considerably: the average solution time is 3.79 times as high when including unavailability ("DynamicRes") as when not including it ("Without unavailability (DynamicRes)"). Furthermore, the solution quality will also be poorer since there are many more timeouts among the instances: when including the unavailability, 11.8\% of the 254 instances time out, whereas only 0.8\% of the instances time out when not including the unavailability.

Excluding the weekends, the machines have between 16 and 18 unavailable days each in 2020. When not including the unavailability in the CG model, the number of fractions scheduled on unavailable days is on average 375.1 per machine during 2020. It could be an option to not include the unavailability in the model, and to instead manually reschedule the fractions for the patients on these days. However, the number of fractions to manually reschedule is very large, and there is no way to ensure that the minimum fractions per week requirements can be fulfilled. 
According to the UK guidelines for the management of treatment interruptions \cite{RCR2019}, there is a wide range of tumor growth rates, and therefore, three different categories with different regulations on maximum number of gaps in the schedule. Category 1 patients have tumor types for which there is evidence that prolongation of treatment affects outcome, and who are being treated radically with curative intent. For these patients, the guidelines state: "Any audit of this category of patient – departmental or national – should show that there was no prolongation of overall treatment time in excess of two days for at least 95\% of the group." The schedules generated by the CG algorithm all fulfill this goal.

 \subsection{Uncertainty in future arrival rates} 
\begin{figure}[ht]
 \center
\includegraphics[scale=0.94]{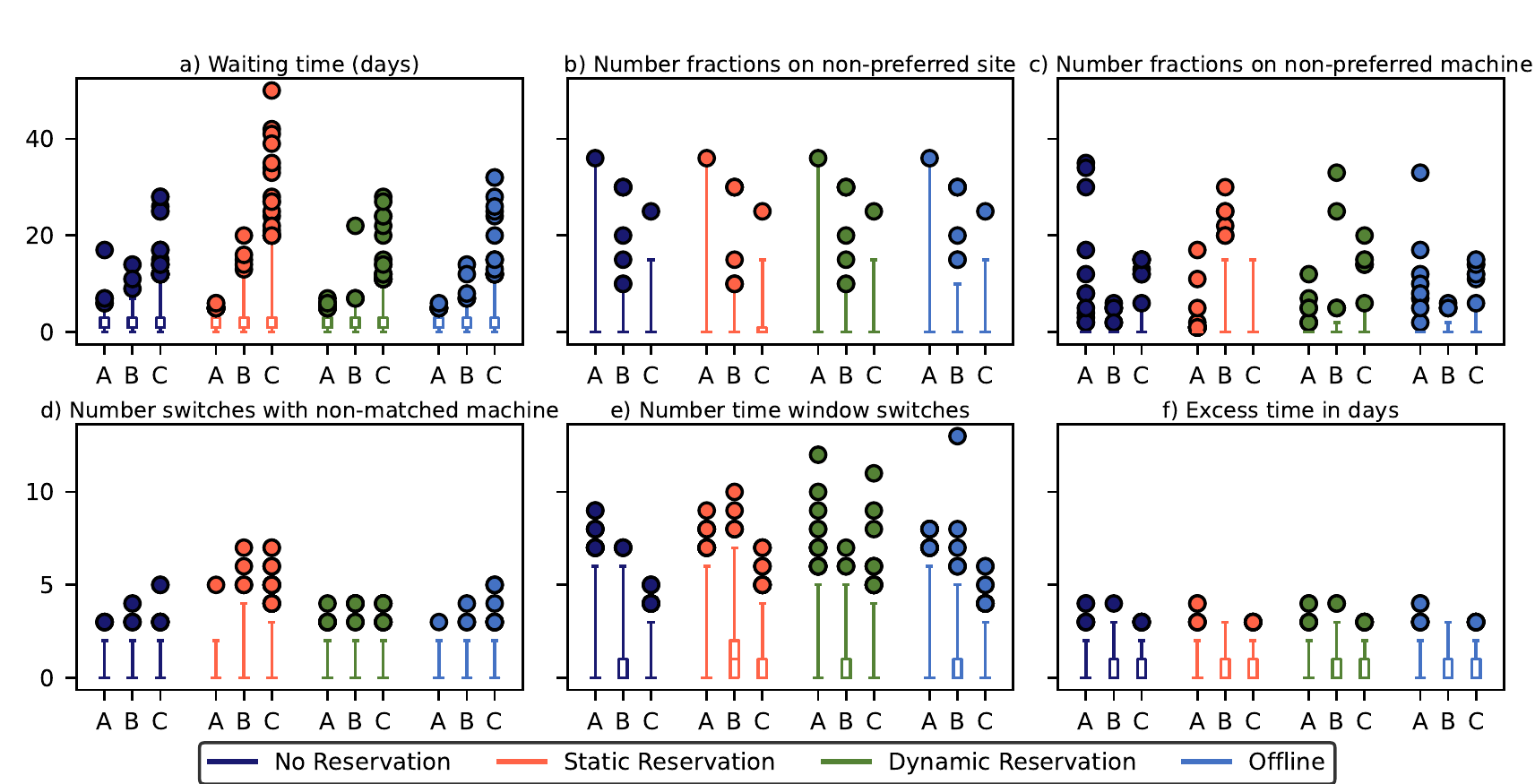}
\caption{Costs relating to the different objectives for priority group A, B and C, with the top 1\% marked as outliers}
\label{Fig:future_all_versions}
\end{figure}

The methods \emph{no reservation}, \emph{static time reservation}, and \emph{dynamic time reservation} are compared to the \emph{offline} solution, which includes the actual future patient arrivals and therefore is the optimal solution. The solution times are presented in Figure~\ref{Fig:with_without_unavailability}. Figure~\ref{Fig:future_all_versions} shows six of the objectives for all versions. 
The solution time is shortest when not taking future patients into account at all (no reservation); this is the only setting in which the CG approach never reaches the time limit. However, in Figure~\ref{Fig:future_all_versions} it is clear that although it is comparable with the other methods in objectives b and d-f, it clearly has the worst waiting time for priority A patients (objective a), and some patients would have to wait almost three weeks for treatment. Furthermore, the priority~A patients also have a higher cost in objective c, i.e., more fractions scheduled on non-preferred machines. 
In contrast, the offline solution naturally has the best performance in the waiting time objective, since it always optimizes the schedules for the actual future events that will occur. The offline solution also demonstrates that there is a trade-off between the objectives; although is has better performance in the waiting time objective than the other methods, this is possible due to more fractions being scheduled on non-preferred machines (objective c) for the priority~A patients.

Figure \ref{Fig:with_without_unavailability} shows that the static time reservation method is much faster than the dynamic time reservation method, and the solution times are comparable to when not reserving time for priority~A patients at all. In Figure~\ref{Fig:future_all_versions}a), it is shown that the static time reservation method has the same performance in waiting time for priority A patients as the offline solution, however, at the cost of priority B and C patients sometimes having to wait unacceptably long for treatment. In objective b-f, Figure~\ref{Fig:future_all_versions} shows that both the static and dynamic time reservation methods have similar performance.

The dynamic time reservation method has the longest solution times as seen in Figure~\ref{Fig:with_without_unavailability}, which could be expected since there is one subproblem for each patient, and therefore the time to solve the subproblems scales somewhat linearly with the number of patients (recall that the dynamic time reservation method adds future patients as placeholder patients that are included in the optimization). The master problem, however, also grows with the number of patients and its increase is sometimes much more than linear due to a combinatoric explosion. Although the solution times are longer and the method often times out before reaching optimality, the dynamic time reservation is the method closest to the offline solution for the waiting times and has similar performance in the other objective functions. This indicates that if the timeout would be longer, the quality of the dynamic time reservation method would likely be even better. Overall, the performance in the objective functions as shown in Figure~\ref{Fig:future_all_versions} demonstrates that the dynamic time reservation method outperforms both the method without reservation and the static time reservation method.

\subsection{Sensitivity analysis for dynamic time reservation method} 
For the dynamic time reservation method, the amount of time reserved (i.e., the number or placeholder priority A patients) is varied to analyze the sensitivity of the solution. The different cases are presented in Table~\ref{Tab:sensitivity_cases}, where the number of future priority~A patients, $\lambda_A$, is varied by $\pm6$ patients per week.  The dynamic time reservation method adds expected future priority A patients for the weeks in the time horizon where the current (actual) patients are assumed to start treatment based on their earliest start day given by $\dMin$. This is usually around 8~weeks, but can also be shorter and longer. Therefore, Table~\ref{Tab:sensitivity_cases} also shows the average total number of future placeholder patients per problem instance. 
\begin{table}[ht]
\caption{The variations in expected future arrivals for sensitivity analysis of the dynamic time reservation method}
\small
\label{Tab:sensitivity_cases}
\setlength\extrarowheight{2pt}
\begin{tabular}{P{4.2cm}P{4.2cm}P{6.8cm}}
\hline
 Number priority A-patients per week & Note & Average total number of future placeholder patients per problem instance\\
\hline
$\lambda_A=30$ & $-16.7\%$ & 240.1 (std: 45.0)\\
$\lambda_A= 36$ & Expected number from data & 288.5 (std: 54.0) \\
$\lambda_A= 42 $&  $+16.7\%$ & 336.7 (std: 63.0)\\
\hline
\end{tabular}
\end{table}

\begin{figure}[ht]
\includegraphics[scale=0.9]{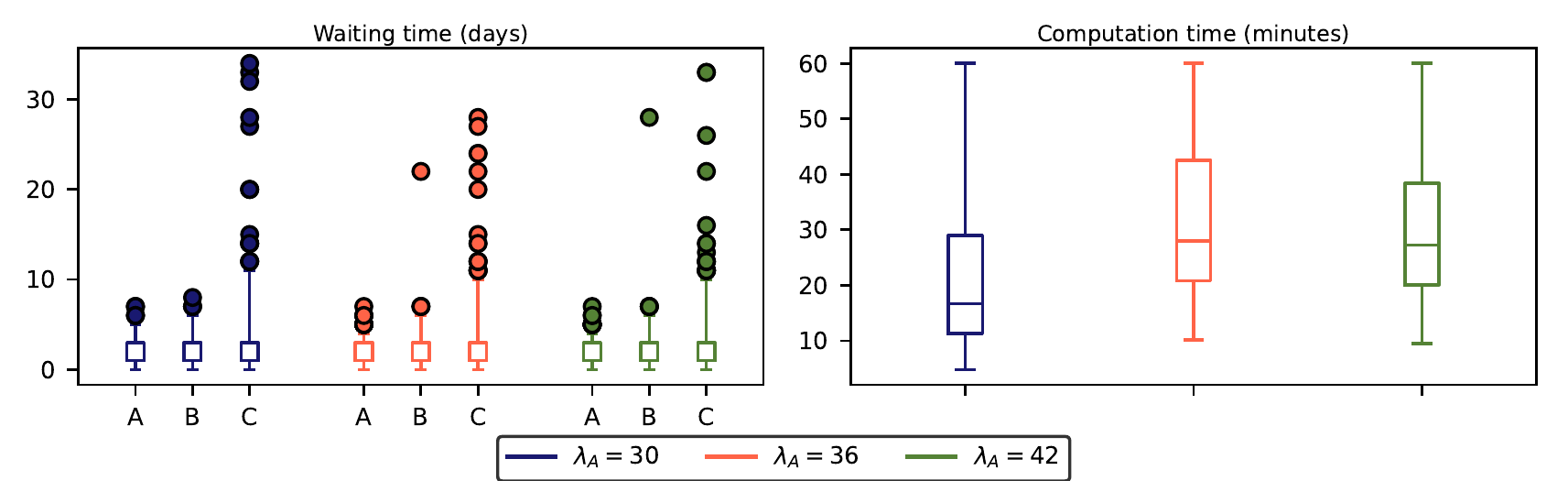}
\caption{Sensitivity analysis for dynamic time reservation, with top 1\% marked as outliers in the waiting time plot}
\label{Fig:sensitivity_dynamic}
\end{figure}

The results from the sensitivity analysis can be seen in Figure~\ref{Fig:sensitivity_dynamic}. The left plot shows the waiting times in days for the different priority groups. It shows that the dynamic time reservation method is robust to fluctuations in the arrival rates, since all three values of $\lambda_A$ perform similarly, with some variations in the top 1\% of the waiting times only. The right plot shows that the computation times for the different values of $\lambda_A$ are also similar, however, it can be seen that for $\lambda_A=42$, the solution times are actually shorter than for $\lambda_A=36$. The most likely explanation to this is that in the CG algorithm, all subproblems (one for each patient) are run only every third iteration, and in the two intermediate iterations only the subproblems that previously had negative reduced costs are run. When $\lambda_A$ is larger, more time is reserved for priority A patients, which means that for the current (actual) patients there is more capacity available in the input schedule, making it more likely that the schedules generated initially by the heuristic can be used in an optimal solution. If the subproblems do not have negative reduced costs, they will not be run in every iteration, saving time considerably.

\subsection{Performance of the CG approach}
As stated in Section~\ref{Sec:ColGen}, the linear relaxation of the master problem is converted to an integer program in the last step, which does not guarantee that the optimal solution is found. To test the performance of the CG approach, an equivalent IP is formulated and run for a subset of ten problem instances for which the CG approach did not time out, since in those cases it is obvious that the solution can be improved. By using the CG solution as a warm start and giving the IP model unlimited runtime, the optimal value to each instance is determined. Table~\ref{Tab:CG_performance} shows the optimality gap between the CG solution and the \textit{proven} optimal solution for a number of different initial columns. It also shows the results if warmstarting the IP formulation from a feasible solution generated by the CG heuristic, for which it can prove optimality in four cases within the 24 hour time limit. The results show that the CG approach performs best for 50 initial columns per patient, but that it is robust to changes in this number. On average, the optimality gap between the solution from the CG approach and the proven optimality is 2.3\% for the 50 column case.

\begin{table}[b]
\small
\caption{The performance of the IP model and the CG approach with different number of initial columns per patient}
\label{Tab:CG_performance}
\setlength\extrarowheight{2pt}
\begin{tabular}{P{0.062\linewidth}P{0.13\linewidth}P{0.13\linewidth}P{0.13\linewidth}P{0.14\linewidth}P{0.25\linewidth}}
\hline
& \multicolumn{5}{P{0.88\linewidth}}{Relative optimality gap between the current best objective value and the \textit{proven} optimal value in percent. In parenthesis is the solution time.}\\
\hline
Instance  &  CG: 30 columns per patient&CG: 50 columns per patient & CG: 75 columns per patient& CG: 100 columns per patient& IP (24h time limit)\\
\hline
46   & 0.5\% (23 min)   & 0.1\% (19 min)   & 0.4\% (20 min)  & 0.0\% (25 min)  &10.3\% (24h)\\
54   & 0.5\% (35 min)   & 0.5\% (23 min)   & 0.3\% (22 min)  & 0.7\% (22 min)  & 10.2\% (24h)\\
74   & 0.0\% (23 min)   & 0.0\% (22 min)   & 0.4\% (25 min)  & 0.0\% (29 min)  & 57.3\% (out of memory at 6.5h)\\
96   & 10.4\% (31 min)& 10.4\% (31 min) & 10.8\% (31 min)& 10.7\% (36 min)& 0.1\% (24h)\\
121 & 0.0\% (19 min)  & 0.0\% (18 min)   & 0.0\% (20 min)   & 0.0\% (20 min)  & 0\% (16.2h)\\
145 & 2.9\% (34 min)  & 1.8\% (33 min)   & 2.5\% (31 min)   & 2.4\% (33 min)  & 0\% (18.8h)\\
165 & 0.0\% (27 min)  & 0.4\% (34 min)   & 0.3\% (24 min)   & 0.4\% (35 min)  & 7.2\% (out of memory at 14.4h)\\
183 & 0.4\% (18 min)  & 0.1\% (22 min)   & 0.1\% (23 min)   & 0.1\% (24 min) & 50.5\% (out of memory at 14.9h)\\
238 & 0.8\% (20 min)  & 0.5\% (25 min)   & 0.5\% (20 min)   & 0.5\% (26 min) & 0\% (1.9h)\\
245 & 9.6\% (34 min)  & 9.3\% (35 min)   & 10.5\% (57 min) & 10.9\% (37 min) & 0\% (42 min)\\
\hline
Average & 2.5\% (26 min) & 2.3\% (26 min) & 2.6\% (27 min) & 2.6\% (29 min) & 13.5\% (872 min)\\
\hline
\end{tabular}
\end{table}

\section{Conclusions}
\label{Sec:conclusions}
It has been shown in numerous studies that long waiting times for RT negatively impacts clinical outcomes. Since the number of linacs in a clinic is limited, the waiting time for treatment is often directly linked to the RT patient scheduling problem, which is currently done manually by the staff in most clinics. 
The main contribution in this paper is that we present an automatic scheduling algorithm for the RTSP based on column generation. The model includes all constraints and objectives necessary for it to work in practice at Iridium Netwerk, a large cancer center with ten linacs in Antwerp, Belgium. 
To the best of our knowledge, this is the first model to consider planned machine unavailability and the constraints and objectives related to the resulting gaps in the schedules. The model is also the first to include specialized treatments, such as consecutive treatments, and non-conventional treatments. The model also supports multiple hospital locations and allows the patients to have preferences on where they want to be treated. To account for uncertainty in the future arrivals of urgent patients, we present a method for dynamic time reservation, and compare it to the static time reservation method that most clinics use today.

The results show that including planned machine unavailability is not straightforward: there are many additional constraints that must be added to ensure that the patients fulfill the minimum fractions per week requirement while the number of gaps in the schedules are minimized. This leads to a large increase in computational time. However, this addition is necessary for the automatic scheduling algorithm to work in practice when scheduling patients on a rolling time horizon. 
Furthermore, when comparing different methods to handle uncertainty in future high-priority patient arrivals, the results show that the dynamic time reservation method outperforms the static time reservation method. This is especially true for lower prioritized patients, as the static method seems to sometimes be too conservative when reserving time for future urgent patients. The dynamic time reservation method frequently reaches the time limit of one hour, thus the quality of the scheduled could possibly be improved further if allowing longer computation times. 
Since the average arrival rates, together with the distribution of the protocols, are taken from the data from 2020, and the dynamic time reservation method is thereafter tested on the 2020 data, there could be a bias and the results could be overvalued. However, the sensitivity analysis show that the dynamic time reservation method is robust to fluctuations in arrival rates, which strengthens the conclusion that this method works very well. Finally, by evaluating the performance of the CG approach compared to an exact method, it can be seen that the CG approach generates schedules that are close to optimal in a reasonable time frame.

\mypar{Future work} The first future step is to compare the automatically generated scheduled with the actual schedules from Iridium Netwerk that were manually constructed. To do this, the manual schedules must be obtained. To make the comparison fair, the same obstacles must be present when using the automatic scheduling algorithm as were in reality, including unplanned unavailability of machines due to failures. 


When accounting for uncertainty in future arrivals, it is possible that the dynamic time reservation method could be improved by not using a static number of placeholder patients to add every week in the time horizon based on average arrivals, but instead using machine learning to predict the future arrivals based on historical arrival rates and the current occupancy in the schedule. To to this, more data is needed, since although a prediction based on the 2020 data can be done, another dataset would be needed for testing. 

To improve the quality of the schedules, it is possible that grouping similar treatments to be scheduled after each other could decrease the times needed for machine setup. Furthermore, the fixation devices used during the RT would not need to be moved between treatment rooms, which could also save time and effort for the medical staff. 

\subsection*{Acknowledgements}
The authors are greatful to Carole Mercier and Geert de Kerf at Iridium Netwerk for valuable insights in the RT scheduling process and help with data gathering. The authors also thank Mats Carlsson at RISE Research Institutes of Sweden for helpful comments about the manuscript.
\clearpage

\begin{appendices}
\section{Pseudocode for schedule construction}
\label{Appendix:heuristic}
The pseudocode for the construction of the initial schedules in the CG algorithm is presented in Algorithm~\ref{Alg:ColGenGreedy}. The function $R: \R^n \to \R^n$ takes a set of elements and returns them in random order. In \eqref{AlgEq:machine_list}, different sets of machines are presented for each patient, that are used in Algorithm~\ref{Alg:ColGenGreedy}. The first set of machines $\mathcal{M}_1$ is the best for the patient, since it adheres to both the site preference and protocol machine preference, and the following are in descending order less well suited for the patient to be scheduled on.  
\begin{equation}
\label{AlgEq:machine_list}
    \begin{aligned}
\mathcal{M}_1^p &=  \SMpref \cap \Mhpref, \qquad &\mathcal{M}_2^p& = \SMpref \cap \Mh\\
 \mathcal{M}_3^p &= \Mhpref, \qquad & \mathcal{M}_4^p &= \Mh
    \end{aligned}
\end{equation}

\begin{algorithm}
\caption{Initial schedule generation procedure for conventional treatments} 
\begin{algorithmic}[1]
\State Initiate $\mathcal{K}_p' = \{\}$ for each $\pinP$
\For {$i=1,\dots,N$} \Comment{$N=50$, number of initial schedules}
\For {$\minM, \dinD, \winW$}
    \State $\StempMDW \gets S_{m,d,w}$ \Comment{Temporary schedule with input schedule as initial state}
\EndFor
    \State $\Hlist \gets \{ R(\mathcal{H}_{\text{priority A}}) \}+ \{R(\mathcal{H}_{\text{priority B}}) \}+ \{R(\mathcal{H}_{\text{priority C}})\}$ \Comment{Protocol list in partly randomized order}
    \For {each protocol $h \in \Hlist$}
        \If {in 50\% of the iterations}  
            \State $\Plist \gets \{R(\mathcal{P}_h)\}$  \Comment{Patient list in randomized order}
        \Else 
            \State $\Plist \gets  \{ p \in \mathcal{P}_h$ sorted in increasing order of $\dMin$\} \Comment{Patient list in deterministic order}
        \EndIf
        \For {each patient $p \in \Plist$}
            \If {in 25\% of the iterations} 
                \State  $\Mlist \gets  \{R(\Mh)\}$ \Comment{Completely randomized}
            \Else
                \State $\Mlist \gets \{R(\mathcal{M}_1^p) \}+ \dots + \{R(\mathcal{M}_4^p)\}$  \Comment{Partly randomized with machine sets from \eqref{AlgEq:machine_list} }
            \EndIf
            \If {$p \in \Ppref$}  \Comment{For patients that have a time window preference}
                \If{in 10\% of the iterations} 
                    \State $\Wlist \gets  \{R(\mathcal{W} \setminus \wPref ) \} + \{\wPref\}$ \Comment{Do not choose preferred time window first}
                \Else
                    \State $\Wlist \gets  \{\wPref\} + \{R(\mathcal{W} \setminus \wPref ) \}$ \Comment{Choose preferred time window first}
                \EndIf
            \Else
                \State $\Wlist \gets  \{ R(\mathcal{W}) \}$
            \EndIf
            \For {each machine $m \in \Mlist$}
                \State  $d_{\textit{start}} \gets$  first available $d$ for protocol $h$ and machine $m$ 
                \State $F_{\textit{scheduled}} \gets 0$ \Comment{Number of fractions scheduled so far \emph{(global variable)}}
                \State $\mathcal{R}_p \gets \{ \}$ \Comment{Used resources for patient \emph{(global variable)}} 
                \For {$d=d_{\textit{start}},\dots,D$} 
                    \State $w_\textit{found} \gets$ False  \Comment{If time window has been found}
                     \If {$d \in \DmU$} \Comment{The day is unavailable for machine $m$}
                         \State $w_\textit{found}  \gets$  \Call{HandleUnavailableDay}{$\Stemp, h, p, m, d, i$}  \Comment{Also updates $\mathcal{R}_p, F_{\textit{scheduled}}$}
                     \Else
                            \For {$w \in \Wlist$}
                                 \If {$dur_p + \StempMDW \leq L_w$} \Comment{There is available capacity}
                                     \State $F_{\textit{scheduled}} \gets F_{\textit{scheduled}}  + 1$
                                     \State  $\mathcal{R}_p \gets \mathcal{R}_p + \{(m, d, w)\}$ \Comment{Save tuple of used resources}
                                     \State $w_\textit{found} \gets$ True
                                     \State \textbf{break} window loop
                                 \EndIf
                            \EndFor
                      \EndIf
                      \If {not $w_\textit{found} $} \Comment{There is no capacity on day $d$, $d+1$ is earliest start day}
                               \State $F_{\textit{scheduled}} \gets 0$ \Comment{Remove fractions}
                               \State $\mathcal{R}_p \gets \{ \}$ \Comment{Reset resources}
                     \EndIf
                  
                  \If {$F_{\textit{scheduled}} = F_p$} \Comment{Found available days for all fractions}
                      \For{$(\hat{m}, \hat{d}, \hat{w}) \in \mathcal{R}_p$}
                          \State $\Stemp_{\hat{m},\hat{d},\hat{w}} \gets \Stemp_{\hat{m},\hat{d},\hat{w}} + \durprest$ \Comment{\emph{Note: Duration is $\durpfirst$ if it is the first fraction}}
                      \EndFor
                      \State $\mathcal{K}_p' \gets \mathcal{K}_p' \cup \{i\}$ \Comment{Extend $\mathcal{K}_p'$ to include new schedule $i$}
                      \State  \textbf{break} day loop and machine loop
                  \EndIf
            \EndFor
            \EndFor
        \EndFor
    \EndFor
\EndFor
\end{algorithmic} 
\label{Alg:ColGenGreedy}
\end{algorithm}

\begin{algorithm}
\caption{Function to handle unavailable day $d$ on machine $m$ for patient $p$ that has protocol $h$}
\label{Alg:UnavailableDayFunction}
\begin{algorithmic}[1]
\Function{HandleUnavailableDay}{$\Stemp, h, p, m, d, i$}
\If {$\fminH  < 5$} \Comment{Patient requires less than 5 fractions per week}
    \If {$i \in \{1,\dots,25\}$} \Comment{In 50\% of cases, try other machine}
        \State \Return \Call{TryBeamMatchedMachine}{$\Stemp,m,d,w$}
    \Else \Comment{In 50\% of cases, try to postpone fraction on same machine}
        \State Get week $g$ that day $d$ is in
        \If {number available days on $m$ in week $g < \fminH$} 
        \State \Return  False  \Comment{Not enough available days on machine $m$ available in week}
        \Else
        \State \Return True \Comment{Postpone fraction to end by stating we have found window but do not update $\mathcal{R}_p$}
        \EndIf
    \EndIf
\Else \Comment{Then patient requires 5 fractions/week}
    \If {$i \in \{1,\dots,25\}$}
        \State $w_{\textit{found}} \gets$  \Call{TryBeamMatchedMachine}{$\Stemp,m,d,w$}
        \If {$w_{\textit{found}}$}
            \State \Return True
         \EndIf
          \LComment{To get 5 fractions this week, only scheduling two on one day remains as an option. First, check in $\mathcal{R}_p$ if two fractions are already scheduled on one day the same week, not allowed to do it twice.}
          \State Get week $g$ that day $d$ is in
          \If {two fractions on some day $\dinD_{g}$ in $\mathcal{R}_p$}
               \State \Return False
          \EndIf
          \LComment{Check if capacity is available on day $d+1$ or day $d-1$. Since two fractions on same day must be scheduled in first and last time window (must be 6 hours apart), on day $d+1$ check capacity in window 1 and on $d-1$ check capacity in window 4. }
          \If {$d+1 \in \mathcal{D}_g$ and $d + 1 \notin \DmU$}
              \If {$dur_p + \Stemp_{m,d,1} \leq L_1$} \Comment{There is available capacity in window 1}
                       \State $F_{\textit{scheduled}} \gets F_{\textit{scheduled}}  + 1$ 
                       \State  $\mathcal{R}_p \gets \mathcal{R}_p + \{(m, d, 1)\}$ \Comment{Save tuple of used resources}
                       \State \Return True 
             \EndIf
          \ElsIf {$d-1 \in \mathcal{D}_g$ and $d - 1 \notin \DmU$}
             \If {$dur_p + \Stemp_{m,d,4} \leq L_4$} 
                       \State $F_{\textit{scheduled}} \gets F_{\textit{scheduled}}  + 1$ 
                       \State  $\mathcal{R}_p \gets \mathcal{R}_p + \{(m, d, 4)\}$ \Comment{Save tuple of used resources}
                       \State \Return True 
             \EndIf
          \EndIf
    \EndIf
\EndIf
\State \Return False
\EndFunction
\end{algorithmic} 
\end{algorithm}

\begin{algorithm}
\caption{Function for trying a beam matched machine}
\label{Alg:BM_machine}
\begin{algorithmic}[1]
\Function{TryBeamMatchedMachine}{$\Stemp,m,d,w$}
\State Get set of beam-matched machines $\binBm$ that $m$ is in
\For {$\hat{m} \in \bm, \hat{m} \neq m$}
    \If {$d \notin \DmUhat$} \Comment{If $d$ not unavailable for machine $\hat{m}$}
        \For {$w \in \Wlist$}
             \If {$dur_p + \Stemp_{\hat{m},d,w} \leq L_w$} \Comment{There is available capacity}
                       \State $F_{\textit{scheduled}} \gets F_{\textit{scheduled}}  + 1$ 
                       \State  $\mathcal{R}_p \gets \mathcal{R}_p + \{(\hat{m}, d, w)\}$ \Comment{Save tuple of used resources}
                       \State \Return True 
             \EndIf
        \EndFor
    \EndIf
\EndFor
\State \Return False
\EndFunction
\end{algorithmic} 
\end{algorithm}

\end{appendices}

\bibliographystyle{unsrt}  
\bibliography{mybibliography}

\end{document}